%% file: FavPoint.tex
\documentclass[reqno]{amsart}

\usepackage{amssymb,latexsym, abbreviations}
\usepackage[mathscr]{eucal}

\usepackage[dvips]{graphics,color}
\usepackage{epsfig}

\usepackage{natbib}

\bibpunct{(}{)}{;}{a}{,}{,}

\newtheorem{lemma}{Lemma}

\theoremstyle{definition}

\theoremstyle{remark}
\newtheorem{remark}{Remark}

\newcommand{\PPP}{\mathbb{P}}
\newcommand{\EEE}{\mathbb{E}}

\newcommand{\nid}{\noindent}
\renewcommand{\qedsymbol}{\ensuremath{\blacksquare}}

\begin{document}
\title[Localization of favorite points]
{Localization of favorite points \\ for diffusion in random
environment}

\author{Dimitrios Cheliotis}

\address{
Department of Mathematics \\ Bahen Center for Information
Technology \\ 40 St. George St., 6th floor \\ Toronto, ON M5S 3G3
\\ Canada  }

\email{dimitris@math.toronto.edu}
\urladdr{http://www.math.toronto.edu/dimitris/}
\thanks{Research partially supported by NSF grant
DMS-0072331  } 

\date{\today}

\begin{abstract}
For a diffusion $X_t$ in a one-dimensional Wiener medium $W$, it
is known that there is a certain process $(b_x(W))_{x\ge0}$ that
depends only on the environment, so that $X_t-b_{\log t}(W)$
converges in distribution as $t\to+\infty$. We prove that, modulo
a small time change, the process $(b_x(W))_{x\ge0}$ is followed
closely by the process $(F_X(e^x))_{x\ge0}$, with $F_X(t)$
denoting the point with the most local time for the diffusion at
time $t$.
\end{abstract}

\maketitle

\section{Introduction}

Consider the diffusion satisfying the formal SDE

\begin{equation}\label{FSDE}
\begin{array}{rl}
dX(t)=&d\omega(t)-\frac{1}{2}W'(X(t))dt, \\
X(0)=&0,
\end{array}
\end{equation}
where $\omega$ is a one-dimensional standard Brownian motion, and
$W$ is a fixed two sided Brownian motion path that we pick before
running the diffusion. $W$ is called the environment. With
probability one, the derivative $W'$ does not exist, but we will
explain later what exactly we mean by such a diffusion.

It in known that there is a real valued process $(b_s(W))_{s>0}$
having great importance for the asymptotic properties of the
diffusion. For example, for almost all $W$, it holds
$(X(t)-b_{\log t}(W))/\log^2 t\to 0$ in probability as
$t\to+\infty$. For stronger results see \cite{HU}, \cite{GO}.

In \cite{DMF}, the authors call this process ``effective
dynamics'' of the motion, and give several properties of the path
$(b_r(W))_{r>0}$. For example, if we define $n(t):=\#$ jumps of
$b$ in the interval $[1,t]$, then $\lim_{n\to+\infty} n(t)/\log
t=4/3$ a.s. And the question is: what does this say about the
diffusion itself? Or, put differently, what object defined in
terms of the diffusion tracks the process $b$? Clearly the
diffusion itself is a much different process than $b$. e.g., it is
recurrent, while $b$ is transient. The next best thing is to find
a process whose value at time $t$ is determined by the knowledge
of $(X_s)_{s\le t}$ only which follows $b$ closely (Note that if
we have the entire path of $X$, then we can completely recover the
process $b$ with probability 1). The one we put forth is the
process of the favorite point of $X$ at time $t$. And the result
justifying this is our theorem, stated below, already announced in
\cite{CH}.

More specifically, to the diffusion $X$ corresponds the local time
process $\{L_X(t,x):t\ge0, x\in\D{R}\}$, which is jointly
continuous and with probability one satisfies

\begin{equation}\label{loctime}
 \int_0^t f(X_s)ds=\int_{\D{R}}f(x)L_X(t,x)dx
\end{equation}
for all $t\ge0$ and any bounded Borel function $f\in\D{R}^\D{R}$.
For a fixed $t>0$, the set $\mathfrak{F}_X(t):=\{x\in\D{R}:
L_X(t,x)=\sup_{y\in\D{R}}L_X (t,y)\}$ of the points with the most
local time at time $t$ is nonempty and compact. Any point there is
called a \textbf{favorite point} of the diffusion at time $t$. One
can prove
that for fixed $t>0$, $\mathfrak{F}_X(t)$ has at most two
elements, and with probability 1, $\mathfrak{F}_X(t)$ has exactly
one element. Also, $Leb(\{t:\mathfrak{F}_X(t) \text{ has two
elements } \})=0$.

Define $F_X:(0,+\infty)\to \D{R}$ with
$F_X(t):=\inf\mathfrak{F}_X(t)$, the smallest favorite point at
time $t$ (what we prove does not change if we define $F_X$ as the
maximum of $\mathfrak{F}_X(t)$).
  Pick any $c>6$, and for any $x$ with
$|x|>1$, define the interval $I(x):=(x-(\log |x|)^c, x+(\log
|x|)^c)$.

Our result says that the processes $F_X(\exp(\cdot))$ and $b$ are
very close. The precise statement is as follows.

\bigskip

\nid \textbf{Theorem.}
\textit{With $\C{P}$-probability 1, there is a $\tau=\tau(W,X)>0$
so that if we label by $(s_n(W, \tau))_{n\ge1}$ the strictly
increasing sequence of the points in $[\tau,+\infty)$ where
$b_{\cdot}$ jumps and $x_n(W, \tau)$ its value in $(s_n,
s_{n+1})$, then there is a strictly increasing sequence $(t_n(X,
W))_{n\ge1}$ converging to infinity so that}

\medskip

(i) $F_X(e^t) \in I(x_n)$ for $n\ge1$, $t\in(t_n, t_{n+1})$.

(ii) $t_n/s_n\to 1$ $\C{P}$-a.s.

\medskip

\nid \textit{We abbreviated $s_n(W, \tau), x_n(W,\tau), t_n(X, W)$
to $s_n, x_n, t_n$.}

\bigskip

The times $\{t_i:i\ge1\}$ will be defined explicitly in the proof
of the theorem. Observe also that for big $x$, the interval $I(x)$
is a relatively small neighborhood of $x$. Thus, the theorem says
that after some point, the function $F_X(\exp(\cdot))$ ``almost
tracks'' the values of the process $b$ with the same order and at
about the same time.

 A corollary of the theorem and of its proof is the
following.

\bigskip

\nid\textbf{Corollary.} \textit{ Let $c>6$ be fixed. With
$\C{P}$-probability 1, there is a strictly increasing map $\gl$
from $[0,+\infty)$ to itself with $\lim_{s\to+\infty} \gl(s)/s=1$
and}
$$|F_X(e^{\gl(s)})-b_s|<(\log s)^c$$
\textit{for all large $s$.}

\bigskip

One can show that $\varlimsup_{s\to+\infty}b_s/s^2\log\log
s=\pi^2/8$ (it follows from the proofs in \cite{HUSHI2}. The
process $b$ is much easier to handle than $X$ or $F_X$.). Using
the corollary, we get
$$\varlimsup_{s\to+\infty}\frac{F_X(e^s)}{s^2\log\log
s}=\frac{8}{\pi^2}.$$ 
\begin{remark}
The diffusion considered above is the continuous time and space
analog of the so-called Sinai's walk, which is a walk taking place
in $\D{Z}$. These two models are connected (see the survey article
\cite{SH1}), and in most respects they behave in exactly the same
way. For Sinai's walk, limiting properties of the process
$(\xi(n))_{n\ge1}$, with $\xi(n)$ being the number of visits paid
to the most visited point by time $n$, have been studied in
\cite{RE} and \cite{DGPS}. More related to our work is
\cite{HUSHI}, where the authors study the process $F^+(n)$ of the
 location of the biggest positive favorite point at time $n$ as
 well as the analog for the diffusion, i.e.,
$F_X^+(t):=\max \mathfrak{F}_X^+(t)$, where
$\mathfrak{F}_X^+(t):=\{x\in[0,+\infty): L_X(t,x)=\sup_{y\ge0}
L_X(t,y)\}$. The results for the diffusion are
$$\varlimsup_{t\to+\infty} \frac{F_X^+(t)}{(\log t)^2 \log \log \log
t}=\frac{8}{\pi^2},$$
and for any non-decreasing function $f>1$,
$$\varliminf_{t\to+\infty} \frac{f(t)}{\log^2 t}F_X^+(t)=\begin{cases} 0 \\ +\infty
\end{cases} \text { a.s } \Leftrightarrow \int^{+\infty} \frac{\log f(t) }{t\sqrt{f(t)} \log t}dt
\begin{cases} = +\infty \\ <+\infty
\end{cases}
$$
The crucial element in the proofs of the above two results is the
fact that $F_X$ and $b$ are closely connected.

Finally, in connection with this paper of Hu and Shi,  we should
mention that in our work we use some of the techniques appearing
there.

\end{remark}

\begin{remark} \label{convinpr}
One can also prove that
$$F_X(t)-b_{\log t} (W)\to 0 \text{ in probability as } t\to+\infty.$$
In the end of this section, we will just give an intuitive
  argument (insufficient however) to show that it is plausible. The statement can be
proved using similar techniques as in the proof of our theorem.
\end{remark}

The paper is organized as follows. In the remaining of this
section, we define some basic objects, and we give an outline of
the proof.
 The main ideas of the paper are contained in Section
 \ref{Thm2sec}, where the theorem is proved assuming that
the environment behaves in the way we expect it to. In Section
\ref{Lemmsec} we show that, with high probability, the environment
indeed behaves the way we assumed.

\bigskip
\subsection{Some definitions}
Now we will define explicitly the diffusion $X$ and the process
$b$. On the space $\C{W}:=C(\D{R})$, consider the topology of
uniform convergence on compact sets, the corresponding $\gs$-field
of the Borel sets, and $\PPP$ the measure on $\C{W}$ under which
the coordinate processes $\{W(t):t\ge0\},\{W(-t):t\ge0\}$ are
independent standard Brownian motions. Also let
$\Omega:=C([0,+\infty))$, and equip it with the $\gs$-field of
Borel sets derived from the topology of uniform convergence on
compact sets. For $W\in\C{W}$, we denote by $\textbf{P}_W$ the
probability measure on $\Omega$ such that $\{X(t):t\ge0\}$, the
coordinate process, is a diffusion with $X(0)=0$ and generator
\[
\frac{1}{2}e^{W(x)}\frac{d}{dx}\left(e^{-W(x)}\frac{d}{dx}\right).
\]
In fact such a diffusion is defined by the formula
$$X_t=A^{-1}(B(T^{-1}(t))),$$
where \begin{align*} A(x)&=\int_0^x e^{W(s)}ds, \\
T(t)&=\int_0^t e^{-2 W (A^{-1}(B(s)))}ds.
\end{align*}
$A$ is the scale function for the diffusion, $B$ is a standard
Brownian motion, and $T$ is a time change. Then consider the space
$\C{W}\times\Omega$, equip it with the product $\gs$-field, and
take the probability measure defined by
\[
d\C{P}(W, X)=d\textbf{P}_W(X)d\PPP(W).
\]
The marginal of $\C{P}$ in $\Omega$ gives a process that is known
as diffusion in a random environment; the environment being the
function $W$.

\nid Throughout the paper, for $x\in\D{R}$, we will use the
notation
\begin{align*}
\tau(x):&=\inf\{t\ge0: X_t=x\},\\
\rho(x):&=\inf\{t\ge0: B_t=x\}
\end{align*}
 for
the hitting times of $X$ and $B$ respectively.

\nid The local time of the diffusion, introduced in
\eqref{loctime}, is given by
\begin{equation} \label{xloctime}
L_X(t,x)=e^{-W(x)}L_B(T^{-1}(t), A(x)),
\end{equation}
where $L_B$ is the local time process of the Brownian motion $B$.
A useful property of $L_B$ is
\begin{equation} \label{locscaling}
\big(L_B(\rho(a),x)\big)_{x\in\D{R}}\overset{law}{=} \big(|a|
L_B(\rho(1), x/a)\big)_{x\in \D{R}},
\end{equation}
where $a\ne0$.
\medskip

\newpage

To define the process $b$, we introduce some terminology. For a
function $f:\mathbb{R}\to\D{R},\ x>0$, and $y_{0}\in\D{R}$, we say
that \textbf{$f$ admits an $x$-minimum at $y_{0}$} if there are
$\ga,\gb\in\mathbb{R}$ with $\ga<y_{0}<\gb$,
$f(y_{0})=\inf\{f(y):y\in[\ga,\gb]\}$ and $f(\ga)\ge f(y_{0})+x$,
$f(\gb)\ge f(y_{0})+x$. We say that \textbf{$f$ admits an
$x$-maximum at $y_{0}$} if $-f$ admits an $x$-minimum at $y_{0}$.
We denote by $R_x(f)$ the set of $x$-extrema of $f$.

It is easy to see that with probability one, for all $x>0$, the
set $R_x(W)$ for a Brownian path has no accumulation point in
$\D{R}$, it is unbounded above and below, and the points of
$x$-maxima and $x$-minima alternate. Thus we can write
$R_x(W)=\{x_k(W,x):k\in\D{Z}\}$, with $(x_k(W,x))_{k\in\D{Z}}$
strictly increasing, and $x_0(W,x)\le 0<x_1(W,x)$. One of
$x_0(W,x), x_1(W,x)$ is a point of local minimum. This we call
$b_{x}(W)$.

For a fixed $x$, the part $W|[x_k, x_{k+2}]$ of the path of $W$
between two consecutive $x$-maxima we call it an
\textbf{$x$-valley}, or simply a valley when the value of $x$ is
understood. The \textbf{depth} of the valley is defined as
$\min\{W(x_k)-W(x_{k+1}), W(x_{k+2})-W(x_{k+1})\}$, and the point
$k_{k+1}$ is called the \textbf{bottom} of the valley.

\subsection{Informal description of the proof}
 Lets look at Figure 1. The diffusion visits first $b_{r^-}$,
then moves on to $b_r$, and then to $b_{r^++}$. The point $\gz_r$
is a point so that $W(\gz_r)-W(b_r)$ is a bit over the depth of
the valley containing $b_r^-$. The point $\eta_r$ is a point so
that $W(\eta_r)-W(b_{r^++})$ is a bit under the depth of the
valley containing $b_r$. It is true that from the time the
diffusion hits $\gz_r$ until the time it hits $\eta_r$ the
favorite point is near $b_r$. This is expected. Lets see an
illustrative calculation that reveals this analytically also. At
time $\tau(\eta_r)$, we compare the local time of any other point
$s$ with the local time of $b_r$. We have
$$
\frac{L_X(\tau(\eta_r),s)}{L_X(\tau(\eta_r),
b_r)}=e^{W(b_r)-W(s)}\frac{L_B(T^{-1}(\tau(\eta_r)),
A(s))}{L_B(T^{-1}(\tau(\eta_r)), A(b_r))}.$$
 But $T^{-1}(\tau(\eta_r))=\rho(A(\eta_r))$. So the above quotient
 equals
$$e^{W(b_r)-W(s)}\frac{L_B(\rho(A(\eta_r)),A(s))}{L_B(\rho(A(\eta_r)),
A(b_r)
)}=e^{W(b_r)-W(s)}\frac{Z_{A(\eta_r)-A(s)}}{Z_{A(\eta_r)-A(b_r)}},$$
where $Z$ is a two dimensional squared Bessel process (Ray-Knight
theorem). We know that for large $t$, we have $Z_t\approx t$. So
that the above ratio is about
\begin{equation} \label{ratio}
e^{W(b_r)-W(s)}\frac{A(\eta_r)-A(s)}{A(\eta_r)-A(b_r)}=\frac{\int_s^{\eta_r}
e^{W(x)-W(s)}dx}{\int_{b_r}^{\eta_r} e^{W(x)-W(b_r)} dx}.
\end{equation}
The dominant contribution to the integrals comes from the points
$x$ where $W(x)-W(s), W(x)-W(b_r)$ are maximum. The exponent in
the integrand of the denominator has maximum $W^{\#}(b_r,
b_{r^++})$. Regarding the numerator, if the point $s$ is in the
valley of $b_r$ and away from $b_r$, the values of $W(x)-W(s)$
will be a bit less than $W^{\#}(b_r, b_{r^++})$ ($W(s)$ will be
larger that $W(b_r)$), while if $s$ is in the valley of
$b_{r^++}$, then by the definition of $\eta_r$, $W(x)-W(s)$ will
always be a bit less than $W^{\#}(b_r, b_{r^++})$. Consequently
the ratio in \eqref{ratio} is less than one for points $s$ that
are reachable by time $\tau(\eta_r)$ and away from $b_r$.

The proof we sketched for the time $\tau(\eta_r)$ is done for all
times in $[\tau(\gz_r), \tau(\eta_r)]$. After the diffusion visits
$\eta_r$, it reaches $\gz_{r^++}$ (see Figure \ref{Fig4}), where
$W(\gz_{r^++})-W(b_{r^++})$ is a little over the depth of the
valley of $b_r$. Then the favorite point will be near $b_{r^++}$.
In this way, we know from about what time the bottom of each
valley starts being the favorite point and when it stops. For each
valley, the scenario we described happens on the complement of a
set with a small probability (which we bound in Lemma
\ref{locprop}).

The next step is to glue together all the time intervals (one
corresponding to each valley). We prove that the probabilities of
all exceptional sets where our scenario fails have finite sum, and
we use the first Borel-Cantelli lemma. The main tool for this step
is Lemma \ref{interpolation} and the bound given in Lemma
\ref{locprop}. Thus we get part $(i)$.

For part $(ii)$, note that the process $F(t)$ jumps from the
vicinity of $b_r$ to the vicinity of $b_{r^++}$ around time
$e^{W^{\#}(b_r, b_{r^++})}$. While the process $b$ jumps exactly
at $W^{\#}(b_r, b_{r^++})$.

\subsection{Comment on Remark \ref{convinpr}} The argument promised goes as
follows. It is known that once the diffusion hits the bottom of a
valley of depth $a$, it takes about time $e^a$ to exit and move on
to another deeper valley. If one starts a diffusion in a valley
and puts reflecting barriers at the ends of the valley, then the
diffusion has invariant measure proportional to $e^{-W(x)}$. For
large times $t$, the local time $L_X(t,x)$ spent at $x$ is
analogous to $e^{-W(x)}$. Now assume we pick a large time $t$. The
diffusion is in a valley with depth at least $\log t$, and the
previous valley had depth $A \log t$ where $A$ is a random
variable with bounded density in $[0,1]$ independent of $t$. Let
$\gep>0$. On a set of probability $>1-c \gep$ we have $A<1-2\gep$.
If we pick an environment $W$ with this property, then on a set of
$\textbf{P}_W$-probability close to one, we have that the
diffusion hits the bottom of the $\log t$-valley before time
$t^{A+\gep}$.
 So the diffusion has spent at
least $t-t^{A+\gep}$ in the current valley. This is much larger
than any period it has spent in any previous single valley
(because $(t-t^{A+\gep})/t^{A+\gep}=t^{1-\gep-A}-1>t^{\gep}$ is
large). In this time, the local times at the sites have come very
close to their limiting values (This holds for the points near the
bottom, the others will be farther from their limiting values, but
whatever local time they collect is significantly less than the
local time corresponding to the points near the bottom). And these
values are maximum near the bottom of the valley.

\section{Proof of the Theorem } \label{Thm2sec}

The strategy of the proof is the following. Let $(y_i)_{i\ge1}$ be
the consecutive values of $b$ in $[1,+\infty)$. First we show that
outside a set $K_i$ of very small probability, the favorite point
$F_X(t)$ is very close to $y_i$ in a time interval $I_i$. This is
the content of lemmata \ref{locprop} and \ref{transition}. Then we
show that the probabilities of the $K_i$'s are summable, and we
use the first Borel-Cantelli lemma. This is accomplished through
Lemma \ref{interpolation} and the bounds given in lemmata
\ref{locprop} and \ref{transition}. Finally, the intervals $I_i$
are such that the right endpoint of $I_i$ coincides with the left
endpoint of $I_{i+1}$

\subsection{A basic lemma.}

\begin{lemma} \label{baselemma}

Assume $W\in C(\D{R})$, and $y_0<0<\gz\le \eta$, $\ga<\gb<\gga$,
$\gz>\gb$, $H>1$, $k_4>9$ are such that

\begin{itemize}

\item[(i)] $P_W(X \text{ hits $y_0$ before } \eta)<H^{-2}.$

\item[(ii)]   $\max_{s\in[y_0, \eta]}W(s)\le H^{3/2}.$

\item[(iii)] $|\eta-y_0|\le H^4.$

\item [(iv)]$\log \frac{A(\eta)-A(\gb)}{A(\gz)-A(\gb)}<H^2$.

\item [(v)]$A(\gb)/A(\gz)>-H^{-5}$.

\item [(vi)] $\sup_{y\in[\gz, \eta]} \sup_{y_0<s<y,\ s\notin(\ga,
\gga)} \int_s^{y} e^{W(t)-W(s)} dt/\int_{\gb}^{y}e^{W(t)-W(\gb)}
dt < H^{-k_4+5}$.

\item [(vii)]$\sup_{y\in[\gz, \eta]} \sup_{y_0<s<y}
(e^{W(s)}\int_\gb^y e^{W(t)-W(\gb)}dt)^{-1} <\exp(H^{3/4}).$
\end{itemize}

\nid Then we have the following quenched probability estimate
$$P_W(\exists t\in [\tau(\gz), \tau(\eta)] \text{ with }
F_X(t)\notin(\ga, \gga))< c H^{-2}.$$
\end{lemma}

\begin{proof}

One setting where this lemma will be applied is shown in Figure
\ref{Fig2}, with $b_r, \gz_r$, $\eta_r$ having the roles of $\gb,
\gz, \eta$ of the lemma.

It is enough to prove that the quantity
$$\sup_{\tau(\gz)\le t\le \tau(\eta)}
\sup_{s\notin(\ga, \gga)}\frac{L_X(t,s)}{L_X(t, \gb)}$$
is greater than one with probability at most $cH^{-2}$. Using
\eqref{xloctime}, we get
$$
\frac{L_X(t,s)}{L_X(t, \gb)}=e^{W(\gb)-W(s)}\frac{L_B(T^{-1}(t),
A(s))}{L_B(T^{-1}(t), A(\gb))}.$$
Note that $T^{-1}(\tau(\gz))=\rho(A(\gz))$ and
$T^{-1}(\tau(\eta))=\rho(A(\eta))$. So we are interested in the
quantity
$$\sup_{\rho(A(\gz))<t<\rho(A(\eta))}\sup_{s\notin(\ga, \gga)}
e^{W(\gb)-W(s)}\frac{L_B(t,A(s))}{L_B(t, A(\gb) )}.$$
\nid Let
$$N:=\min\{ k\in \D{Z}: k\ge (\log 2)^{-1}\log(\frac{A(\eta)-A(\gb)}{A(\gz)-A(\gb)})\},$$
$u_k=2^k (A(\gz)-A(\gb))+A(\gb) $ for $k=0,1,...,N-1$, and
$u_N=A(\eta)$.

\nid There are unique points $\gz=:p_0<p_1<\ldots<p_N:=\eta$ such
that $u_k=A(p_k).$

\nid Let $A_1:=[\ X \text{ hits $\eta$ before } y_0\ ]$. Then
\begin{multline*}
P_W(\sup_{\rho(u)<t<\rho(v)}\sup_{s\notin(\ga, \gga)}
e^{W(\gb)-W(s)}\frac{L_B(t,A(s))}{L_B(t, A(\gb) )}>1)\\
\le P_W(A_1^c)+H^2 \sup_{k<N}
P_W(\{\sup_{\rho(u_k)<t<\rho(u_{k+1})}\sup_{s\notin(\ga, \gga)}
e^{W(\gb)-W(s)}\frac{L_B(t,A(s))}{L_B(t, A(\gb))}>1\}\cap A_1).
\end{multline*}
Now
\begin{multline}\label{eqsplit1}
P_W(\{\sup_{\rho(u_k)<t<\rho(u_{k+1})}\sup_{s\notin(\ga, \gga)}
e^{W(\gb)-W(s)}\frac{L_B(t,A(s))}{L_B(t, A(\gb))}>1\}\cap A_1)
 \\
\le P_W(\{\sup_{s\notin(\ga,
\gga)}\frac{e^{W(\gb)-W(s)}}{u_{k+1}-A(\gb)}L_B(\rho(u_{k+1}),A(s))>
\frac{1}{H^4}\}\cap A_1)+P_W(\frac{1}{u_{k+1}-A(\gb)}
L_B(\rho(u_k), A(\gb))<\frac{1}{H^4}).
\end{multline}
To bound the first term, we observe that the local time appearing
in the expression is zero for $s\ge p_{k+1}$, and we use the
Ray-Knight theorem to get
\begin{multline*}
\{\frac{e^{W(\gb)-W(s)}}{u_{k+1}-A(\gb)}L_B(\rho(u_{k+1}),A(s)):s<
p_{k+1}\}\overset{law}{=}
\{\frac{e^{W(\gb)-W(s)}}{u_{k+1}-A(\gb)}\tilde Z_{u_{k+1}-A(s)}:s<
p_{k+1}\}\\ =\{(\int_{\gb}^{p_{k+1}} e^{W(y)-W(\gb)}dy)^{-1}
e^{-W(s)}\tilde Z_{\int_s^{p_{k+1}} e^{W(y)}dy}: s<p_{k+1}\},
\end{multline*}
%
%
where $(\tilde Z_s)_{s\ge0}$ is a two dimensional squared Bessel
process up to time $u_{k+1}$ and then zero dimensional squared
Bessel process. Let also $Z$ be the two dimensional squared Bessel
process which is run with the same Brownian motion as $\tilde Z$.
Then with probability one,
\begin{equation}\label{domination}
 Z_t\ge \tilde Z_t \text{ for all } t\ge0
\end{equation}
by Theorem 3.7 of Chapter IX in \cite{RY}. The function $\rho$
required by that theorem is in our case $\rho(x)=x$ for all
$x\ge0$. So that
$$\{\frac{e^{W(\gb)-W(s)}}{u_{k+1}-A(\gb)}L_B(\rho(u_{k+1}),
\frac{A(s)}{u_k}):s\le p_{k+1}\}\overset{law}{\le} \{
(\int_{\gb}^{p_{k+1}} e^{W(y)-W(\gb)}dy)^{-1} e^{-W(s)}
Z_{\int_s^{p_{k+1}} e^{W(y)}dy}:s\le p_{k+1}\}.$$
Let $$x_1(s):=e^{W(\gb)-W(s)} ,\ x_2(s):=\int_s^{p_{k+1}}
e^{W(y)}dy=A(p_{k+1})-A(s).$$
We know that \begin{equation} \label{envconditions}
\sup_{y_0<s<p_{k+1}, s\notin(\ga, \gga)} \frac{x_1(s)
x_2(s)}{x_2(\gb)}\le H^{-k_4+5} \text{ and } \sup_{y_0<s<p_{k+1},
s\notin(\ga, \gga)} \frac{x_1(s)}{x_2(\gb)}\le \exp(H^{3/4}).
\end{equation}
The first quantity in \eqref{eqsplit1} is bounded by
\begin{multline}
P_W(\{\sup_{s\notin(\ga, \gga),\
s<p_{k+1}}\frac{x_1(s)}{x_2(\gb)}Z_{x_2(s)}>\frac{1}{H^4}\}\cap
A_1 )
\\\le  P_W(\{\sup_{\substack{s\notin(\ga, \gga) \\ x_2(s)\le
e^{-H^2},\  s<p_{k+1}}} \frac{x_1(s)}{x_2(\gb)}Z_{x_2(s)}
>\frac{1}{H^4} \} \cap
A_1 )  + P_W(\{\sup_{\substack{s\notin(\ga, \gga) \\
x_2(s)\ge e^{-H^2}}}
\frac{x_1(s)x_2(s)}{x_2(\gb)}\frac{Z_{x_2(s)}}{x_2(s)}>\frac{1}{H^4}
\} \cap
A_1) \\
\le
 P_W(\sup_{t\le  e^{-H^2}
}e^{H^{3/4}}Z_t>\frac{1}{H^4}) + P_W(\sup_{ e^{-H^2}\le
t \le H^4 e^{H^2} }\frac{Z_t}{t}> H^{k_4-9})\\
\label{estimate1} =
 P_W(\sup_{t\le 1
}Z_t>\frac{e^{H^2-H^{3/4}}}{H^4}) + P_W(\sup_{ 1\le t \le H^4
e^{2H^2} }\frac{Z_t}{t}> H^{k_4-9})
\end{multline}
To justify the last inequality, we use the bound $x_2(s)\le H^4
e^{H^2}$ and \eqref{envconditions}. By well known property of
Brownian motion, $\PPP(\sup_{t\le 1}Z_t>a)\le
4\PPP(B_1>\sqrt{\frac{a}{2}})<4 e^{-a/4}.$ So the first term in
\eqref{estimate1} is bounded by
$4\exp(-\frac{e^{H^2-H^{3/4}}}{4H^4}).$ To bound the last term in
\eqref{estimate1}, we use Lemma \ref{LILbound} with the choices
$a=\sqrt{\frac{H^{k_4-9}}{2}}, M=H^4e^{2H^2}$, and $c=2$, and we
find that it is bounded by $ \exp\left({-H^{k_4-9}/17}\right)$ for
large $H$.

\nid The second term in \eqref{eqsplit1} is bounded as follows. If
$\gb\ge 0$, then we use the Ray-Knight theorem to obtain
$$P_W(\frac{1}{u_{k+1}-A(\gb)} Z_{u_k-A(\gb)}<\frac{1}{H^4})=
P_W(\frac{u_k-A(\gb)}{u_{k+1}-A(\gb)} Z_1<\frac{1}{H^4})\le P_W(
Z_1<\frac{2}{H^4})\le \frac{1}{H^4}.$$
Here $Z$ is a two dimensional squared Bessel process. And we used
also the scaling property of $Z$, the fact that $Z_1$ has a
density bounded by 1/2 (it is exponential with mean 2), and that
$(u_{k+1}-A(\gb))/(u_k-A(\gb))\le 2$. For the case $\gb<0$, we
will need the inequality $(u_{k+1}-A(\gb))/u_k<4$. This translates
to $2^{k+1}(A(\gz)-A(\gb))+4A(\gb)>0$. The last quantity is enough
to be positive for $k=0$. Then the inequality becomes
$A(\gz)+A(\gb)>0$ which holds because of (v). Thus using
\eqref{locscaling} and $(u_{k+1}-A(\gb))/u_k<4$, we get
$$\frac{1}{u_{k+1}-A(\gb)} L_B(\rho(u_k),
A(\gb))\overset{law}{=}\frac{u_k}{u_{k+1}-A(\gb)} L_B(\rho(1),
\frac{A(\gb)}{u_k})>\frac{1}{4} L_B(\rho(1),
\frac{A(\gb)}{u_k}).$$
So that the term is bounded by $P_W( L_B(\rho(1),
\frac{A(\gb)}{u_k})<4/H^4)$. The process $\tilde Z_s:=L_B(\rho(1),
1-s), s\ge0,$ is up to time 1 a two dimensional square Bessel
process, and after that a zero dimensional square Bessel process.
Let $(\hat Z)_{s\ge0}$ be a two dimensional square Bessel process.
Then the comparison Theorem IX.3.7 in \cite{RY}, the fact that
$\tilde Z_1, \hat Z_1$ are exponential with mean 2, and the
assumption $0>A(\gb)/A(\gz)>-H^{-5}$ give
\begin{multline*}P_W(\tilde Z_{1-A(\gb)/u_k}<\frac{4}{H^4})\le  P_W(\tilde
Z_1<\frac{8}{H^4})+P_W(|\tilde Z_1-\tilde
Z_{1-A(\gb)/u_k}|>\frac{4}{H^4}) \\ \le \frac{4}{H^4}+P_W(\hat
Z_{-A(\gb)/u_k}>\frac{4}{H^4})< \frac{4}{H^4}+P_W(\hat
Z_1>-\frac{4}{H^4} \frac{u_k}{A(\gb)})=
\frac{4}{H^4}+\frac{1}{2}\exp(\frac{2}{H^4} \frac{u_k}{A(\gb)})<c
H^{-4}.
\end{multline*}

\nid Putting all estimates together, we get the bound
$$P_W(A_1^c)+H^2\Big(C\ H^{-4}+4\exp(-\frac{e^{H^2-H^{3/2}}}{4H^5})
+\exp\left({-H^{k_4-9}/17}\right)\Big)
 <CH^{-2}.$$
 We used the fact that $k_4>9$.
\end{proof}

\subsection{Two lemmata}

Let $r>0$ be fixed, and $x_0(W,r), x_1(W,r)$ be the $r$-extrema
around zero, as in the definition of $b$. Assume that
$b_r=x_1(W,r)$. Otherwise, all the definitions following should be
applied to the path $(s\mapsto W(-s))$. Let $r^-=\sup\{x<r: b_x\ne
b_r\}, r^+=\inf\{x>r: b_x\ne b_r\}$ the points where $b$ jumps
just before and after $r$ respectively. Note that the function $
(r\mapsto b_r)$ (and similarly all functions $(r\mapsto x_i(W,r)),
i\in \D{Z}$) is step and left continuous. So for the next value of
it after $r$ we will use the notation $b_{r^++}$, i.e., the right
limit of $b$ at the point $r^+$. The probability that $r^-=r$ or
$r^+=r$ is zero, so in the following we assume that $r^-<r<r^+$.

\nid Since $b_r>0$, it holds $b_{r^-}<b_r$. For $x, y\in\D{R}$, we
define
$$
W^{\#}(x,y)=\sup \{W(s)-W(t): (t-s)(x-y)\ge0\}.
$$
This is connected with the time it takes for the diffusion
starting at $x$ to reach $y$.

\nid In the following, we will use three constants $k_1, k_2,
k_3$. Our assumption for them is that $k_1, k_3>9, k_2\ge 18$. We
prefer not to choose values for them so that their role in the
proof is clearer.

\nid Let

\begin{align}
j_r&:=\sup\{s<b_r:W(s)-W(b_r)=r\},\\
l_r&:=\inf\{s>b_r:W(s)-W(b_r)=r\}.
\end{align}
and

\begin{align}
\ga_r&:=\inf\{t>j_r:W(t)-W(b_r)<k_1\log r\},\\
\gga_r&:=\sup\{t<l_r:W(t)-W(b_r)<k_1\log r\}.
\end{align}
So that, for $t\in [j_r,l_r]\setminus (\ga_r, \gga_r)$, we have
\begin{equation}\label{separation}
 W(t)\ge W(b_r)+k_1\log r.
\end{equation}
Also define

\begin{align*}
\gz_r:&= \inf\{t>b_r:W(t)-W(b_r)\ge W^\#(b_{r^-}, b_r)+k_2\log
W^\#(b_{r^-}, b_r) \},\\
\eta_r:&=
\begin{cases} \ \inf\{t>b_{r^++}:W(t)-W(b_{r^++})\ge W^\#(b_r, b_{r^++})-k_3\log
W^\#(b_r, b_{r^++})\}
 & \text{ if } b_{r^++}>0,
\\ \sup\{t<b_{r^++}:W(t)-W(b_{r^++})\ge
W^\#(b_r, b_{r^++})-k_3\log W^\#(b_r, b_{r^++}\} & \text{ if }
b_{r^++}<0,
\end{cases}
\end{align*}
and for any $r\in\D{R}$, let $\tau(r):= \inf\{t>0:X_t=r\}.$

The two cases $b_{r^++}>0, b_{r^++}<0$ are shown in Figures
\ref{Fig2}, \ref{Fig3} respectively along with other points which
are introduced in the proof of Lemma \ref{locprop}.

\begin{figure}[ht]
\begin{center}
\resizebox{9cm}{!}{\input{FavP2.pstex_t}}\caption{The scenario
$b_r>0, b_{r^++}>0$} \label{Fig2}
\end{center}
\end{figure}

\begin{figure}[ht]
\begin{center}
\resizebox{9cm}{!}{\input{FavP3.pstex_t}}\caption{The scenario
$b_r>0, b_{r^++}<0$} \label{Fig3}
\end{center}
\end{figure}

\medskip

\nid It will be shown that with high probability $\tau(\gz_r)<
\tau(\eta_r)$, and the main claim is that from time $\tau(\gz_r)$
to $\tau(\eta_r)$, the favorite point is around $b_r$. The precise
statement is the following lemma.

\begin{lemma} \label{locprop}
For $r>0$,
\begin{equation}
\C{P}\Big(\tau(\gz_r)>\tau(\eta_r) \text{ or there is a $t\in
[\tau(\gz_r), \tau(\eta_r)]$ such that } F(t)\notin(\ga_r,
\gamma_r)\Big)<c r^{-1/4}.
\end{equation}
\end{lemma}

\begin{proof}

We are always under the assumption $b_r>0$. We consider two cases.

\medskip

\nid \textsc{Case} 1: $b_{r^++}>0.$

\medskip

We will apply Lemma \ref{baselemma} for the path $W$ and the
choice $$(y_0, \ga, \gb, \gga, \gz, \eta, H):=(x_0(W,r^++), \ga_r,
b_r, \gga_r, \gz_r, \eta_r, r).$$
Call $A$ the event in the probability of the statement, and $E_1$
the event that one of (1)-(7) fails. The probability of $A$ is
bounded by
$$\PPP(E_1)+\C{P}(A \cap E_1^c) \le cr^{-1/4}.
$$
The first bound is proved in Lemma \ref{thm2lemma}, the second
follows by Lemma \ref{baselemma}.

\medskip

\nid \textsc{Case} 2: $b_{r^++}<0$.

\nid Let
\begin{align}
\tilde \gz_r:&=\sup\{s<0: W(s)-W(b_r)\ge W^{\#}(b_r,0)+2k_2\log W^{\#}(b_r,0)\}, \\
\hat \gz_r:&=\inf\{s>\gz_r: W(s)-W(b_r)\ge W^{\#}(b_r,0)+3k_2\log
W^{\#}(b_r,0) \},
\end{align}
and
$$A_3:=[
 \tau(\gz_r)<\tau(\tilde
\gz_r)<\tau(\hat \gz_r)<\tau(\eta_r)].$$ See Figure \ref{Fig3}.
Then
$$ \C{P}(\tau(\gz_r)> \tau(\eta_r) \text{\  or } F([\tau(\gz_r),
\tau(\eta_r)]) \nsubseteq (\ga_r, \gamma_r) ) \le
\C{P}(\{F([\tau(\gz_r), \tau(\eta_r)]) \nsubseteq (\ga_r,
\gamma_r)\} \cap A_3)+\C{P}(A_3^c).
$$
As Lemma \ref{hitlemma} shows, $\C{P}(A_3^c)<C r^{-2}.$ The first
quantity is bounded by
$$\C{P}(F([\tau(\gz_r), \tau(\hat \gz_r)]) \nsubseteq (\ga_r, \gamma_r))+
\C{P}( F([\tau(\tilde \gz_r), \tau(\eta_r)])\nsubseteq (\ga_r,
\gamma_r)).$$
Both of these two probabilities are bounded with the use of Lemma
\ref{baselemma}.

\nid For the first, we apply Lemma \ref{baselemma} for the path
$W$ and the choice $$(y_0, \ga, \gb, \gga, \gz, \eta,
H):=(x_0(W,r), \ga_r, b_r, \gga_r, \gz_r, \hat \gz_r , r).$$
Working as in Case 1, we obtain the required bound.

 \nid For the second, we apply Lemma
\ref{baselemma} for the path $W^*:=W(-\cdot)$ and the choice
$$(y_0, \ga, \gb, \gga, \gz, \eta, H):=(-x_1(W,r^++), -\gga_r,
-b_r, -\ga_r, -\tilde \gz_r, -\eta_r, r).$$
We will use the notation $X^W$ for the diffusion run in the fixed
environment $W$. Let $E_2$ be the event that, with these choices,
one of (i)-(vii) fails. As in Case 1, we use the bound on
$\PPP(E_2)$ given in Lemma \ref{thm2lemma} to get that for $W$
outside $E_2$, we have
$$P_W(F_{X^{W^*}}([\tau(-\tilde \gz_r), \tau(-\eta_r)])\nsubseteq (-\gga_r,
-\ga_r))<cr^{-1/4}$$
But $X^{W^*}\overset{law}{=} -X^W$. So that
$$P_W (F_{X^W}([\tau(\tilde \gz_r), \tau(\eta_r)])\nsubseteq (\ga_r,
\gga_r))<cr^{-1/4}$$
as required.
\end{proof}

In the time interval $[\tau(\eta_r), \tau(\gz_{r^++})]$, we will
show that $F$ jumps from a neighborhood of $b_r$ to a neighborhood
of $b_{r^++}$. That is, if we let
$$K_r:=\left\{
\begin{array}{c} \text{
There is a time $z_r\in[\tau(\eta_r), \tau(\gz_{r^++})]$ so that}
 \\
 F(t)\in (\ga_r, \gga_r)$ \text{ for } $t\in[\tau(\eta_r) ,
z_r),\\
\text{ and $F(t)\in (\ga_{r^++}, \gga_{r^++})$ for $t\in (z_r,
\tau(\gz_{r^++})]. $}
\end{array}\right\},$$
then the following holds.
\begin{lemma} \label{transition}
For $r>0$,
$$\C{P}(K_r^c)<c r^{-1/4}.$$
\end{lemma}

\begin{proof}

Assume that $b_{r^++}>0$. Then $b_r<b_{r^++}<\eta_r<\gz_{r^++}$.
Let
\begin{align*}
m_r:&=\sup\{s<b_{r^+}:W(s)-W(b_r)=W^{\#}(b_r, b_{r^++})\},\\
z:&=\sup\{s<b_{r^++}:W(s)-W(b_{r^+})\ge W^{\#}(b_r, b_{r^++})+k_2
\log W^{\#}(b_r, b_{r^++})+k_2 \log r\},
\end{align*}
\begin{figure}[ht]
\begin{center}
\resizebox{8cm}{!}{\input{FavP4.pstex_t}}\caption{} \label{Fig4}
\end{center}
\end{figure}
and remember that $W(\gz_{r^+})-W(b_{r^+})=W^{\#}(b_r,
b_{r^++})+k_2 \log W^{\#}(b_r, b_{r^++})$.

\nid Also call $Y_t:=X_{\tau(z)+t}$ the diffusion after the time
$\tau(z)$. All objects defined for $X$ (e.g., the local time, the
process of the favorite point) are defined analogously for $Y$.

\nid Define the events
\begin{align*}
\gS_0&:=[z>m_r],\\
\gS_1&:=[F_X(\tau(\eta_r))\in (\ga_r, \gga_r)], \\
\gS_2&:=[F_X(\tau(\gz_{r^++}))\in (\ga_{r^++}, \gga_{r^++})],\\
\gS_3&:=[(X_{\tau(\eta_r)+s})_{s\ge0} \text{ hits } \gz_{r^++}
\text{ before } z], \\
\gS_4&:=[ F_Y([\tau_Y(\eta_r), \tau_Y(\gz_{r^++})]) \subset
(\ga_{r^++}, \gga_{r^++})].
\end{align*}
On $\gS_1\cap \gS_3 \cap \gS_4$ we claim that
\begin{equation} \label{twocases}
t\in [\tau(\eta_r), \tau(\gz_{r^++})] \Rightarrow F_X(t)\in
(\ga_r, \gga_r)\cup (\ga_{r^++}, \gga_{r^++}).
\end{equation}
Let $t\in [\tau(\eta_r), \tau(\gz_{r^++})]$. Points in $(-\infty,
z]$ collect local time only from the part $(X_s)_{s\le
\tau(\eta_r)}$ of the path $(X_s)_{s\le \tau(\gz_{r^++})}$ by the
definition of $\gS_3$. And by the definition of $\gS_1$, the ones
with the most local time are in $(\ga_r, \gga_r)$.  Points in $[z,
+\infty)$ collect local time only from the part $(X_s)_{\tau(z)\le
s\le t}$ of the path. And by the definition of $\gS_4$, we know
that out of them, the ones with the most local time at time $t$
are in $(\ga_{r^++}, \gga_{r^++})$. This proves our claim.

\nid On $\gS_0\cap\gS_1\cap \gS_2 \cap \gS_3 \cap \gS_4$ we know
that $F(\tau(\eta_r))\in (\ga_r, \gga_r)$, $F(\tau(\gz_{r^++}))\in
(\ga_{r^++}, \gga_{r^++})$, and from time $\tau(\eta_r)$ to
$\tau(\gz_{r^++}), X$ does not visit $(\ga_r, \gga_r)$. These
combined with \eqref{twocases} show that $\gS_0\cap\gS_1\cap \gS_2
\cap \gS_3 \cap \gS_4\subset K_r$.

The proof will be completed after we bound the probability of
$(\gS_0\cap\gS_1\cap \gS_2 \cap \gS_3 \cap \gS_4)^c$.

\nid Lemmata \ref{locationslemma}.3 and \ref{locprop} give the
bound for $\C{P}(\gS_0^c)$ and $\C{P}(\gS_1^c)$ respectively. To
bound $\C{P}(\gS_2^c)$, we apply Lemma \ref{baselemma} with the
choice
$$(y_0, \ga, \gb,
\gga, \gz, \eta, H):=(x_0(W, r^++), \ga_{r^++}, b_{r^++},
\gga_{r^++}, \gz_{r^++}, \gz_{r^++}, r).$$
Lemma \ref{hitlemma} shows that $\C{P}(\gS_3^c)<r^{-2}$.

\nid For $\gS_4^c$, we write $\gS_4^c \subset (\gS_4^c \cap
\gS_0)\cup \gS_0^c$  The probability of $\gS_4^c \cap \gS_0$ is
bounded with the use of Lemma \ref{baselemma} for the environment
$W^z:=W(z+\cdot)$ and with the choice
$$(y_0, \ga, \gb,
\gga, \gz, \eta, H):=(m_r-z, \ga_{r^++}-z, b_{r^++}-z,
\gga_{r^++}-z, \eta_r-z, \gz_{r^++}-z, r).$$
Let $E_4$ be the event that, with this choice, one of (i)-(vii)
fails. Lemma \ref{thm2lemma} shows that $\PPP(E_4)<c r^{-1/4}$,
and as in Lemma \ref{locprop}, we show that $\C{P}(\gS_4^c \cap
\gS_0)<c r^{-1/4}.$ This finishes the proof.
\end{proof}

\subsection{Proof of the main results}

\subsubsection{Proof of the Theorem:} Using lemmata
\ref{locprop} and \ref{transition}, we now prove our theorem. Pick
any $a\in(0,1/2)$ and let $r_k=\exp(k^a)$ for $k\ge1$. Then
$\sum_{k=0}^{+\infty} (r_k)^{-1/4}<+\infty$, and lemmata
\ref{locprop} and \ref{transition} imply that there is a $k_0>0$
so that for $k\ge k_0$
\begin{equation}\label{local} F_X(t)\in (\ga_{r_k}, \gga_{r_k}) \text{ for
all } t\in (z_{r_k^-}, z_{r_k}).
\end{equation}
Indeed, this holds for $t\in [\tau(\gz_{r_k}), z_{r_k})$ clearly.
Now let $j=\max\{n:r_n< r_k^-\}$. Then $F_X(t)\in (\ga_{r_j^+},
\gga_{r_j^+}) \text{ for all } t\in (z_{r_j}, \tau(\gz_{r_j^+})]$.
But $r_j^+=r_k^-$, so $\gz_{r_j^+}=\gz_{r_k}$,
$z_{r_j}=z_{r_k^-}$, and $(\ga_{r_j^+}, \gga_{r_j^+}) \subset
(\ga_{r_k}, \gga_{r_k})$.

Let $(s_n(W))_{n\ge 1}$ be the increasing sequence of the point
where $b|[r_{k_0}, +\infty)$ jumps, $x_n:=b_{s_{n+1}}$, the value
of $b$ in $(s_n(W), s_{n+1}(W)]$, and $t_n:=\log z_{s_n}$. Using
Lemma \ref{interpolation}, we may assume that between any two
terms from $(s_n(W))_{n\ge 1}$ there is a term from $(r_k)_{k\ge
k_0}$. So that, for $n\ge 1$ there is a $k$ with $s_n<r_k\le
s_{n+1}$. Then $z_{r_k^-}=z_{s_n}=e^{t_n},
z_{r_k}=z_{s_{n+1}}=e^{t_{n+1}}$, and $F_X((e^{t_n},
e^{t_{n+1}}))\subset (\ga_{r_k}, \gga_{r_k})$ because of
\eqref{local}. Also $(\ga_{r_k}, \gga_{r_k})\subset I(b_{r_k})$
because of Lemma \ref{intervlemma}. This proves (i) of the
theorem.

For the second claim of the theorem, observe that $t_n/s_n=\log
z_{s_n}/s_n$, and $z_{s_n}\in[\tau(\eta_{s_n}),
\tau(\gz_{s_n+})]$. One can see that
$$\varlimsup_{n\to+\infty} \frac{\log \tau(\gz_{s_n+})}{s_n}
\le 1 \text{ and } \varliminf_{n\to+\infty} \frac{\log
\tau(\eta_{s_n+})}{s_n} \ge 1.$$
The proof of these two is done by modifying the proof of (4.7),
(4.8), (4.11) in \cite{HUSHI2}. Since no new idea is involved, we
omit it.
 \hfill \qedsymbol

\subsubsection{Proof of the Corollary:}

Let $\gl(s)=(t_1/s_1) s$ for $s\in[0,s_1]$. Then in any interval
$[s_n, s_{n+1}]$ with $n\ge1$, $\gl(s)$ is defined as the unique
increasing map of the form $\gga_n s+\gd_n$ mapping that interval
to $[t_n, t_{n+1}]$. Since on $[s_n, s_{n+1}]$ the function
$\gl(s)/s=a_n+\gd_n s^{-1}$ is monotone, it maps $[s_n, s_{n+1}]$
to the interval with endpoints $t_n/s_n, t_{n+1}/s_{n+1}$. It
follows that $\lim_{s\to+\infty}\gl(s)/s=1$. Then for all large
$s$, the theorem says that $|F_X(e^{\gl(s)})-b_s|<(\log |b_s|)^c$.
It is easy to prove that, with probability 1, $\log|b_s|<3\log s$
for all large $s$. [This is similar with the proof the proof of
\eqref{lowerb}. We show $\log \gb_s^+<3\log s$ for large $s$ (see
next section for notation). The basic ingredient is that for large
$A>0$, it holds $\PPP(\gb_s^+>s^2 A)=\PPP(\gb_1^+>A)\le C
e^{-A/2}$. This follows from 5 of Lemma \ref{tailslemma} and the
fact that $b_s\in\{\gb_s^+, \gb_s^-\}$. We don't care for the best
bound $Ce^{-A\pi^2/8}$.] This finishes the proof.
\hfill \qedsymbol

\section{Some Lemmata.} \label{Lemmsec}
In this section, we prove several facts we needed in the proof of
Lemma \ref{locprop} and of the theorem.

\nid First we show an alternative way of computing the process
$b$. Call $W^+$ the process $(W(s):s\ge0)$ and $W^-$ the process
$(W(-s):s\ge 0)$. For $r>0$, let
\begin{align*}
\tau_r^+&:=\min\{s\ge0:W^+(s)-\ul{W}^+(s)=r\},\\
\gb_r^+&:=\min\{s\ge0:W^+(s)=\ul{W}^+(\tau_r^+)\},\\
\tau_r^-&:=\min\{s\ge0:W^-(s)-\ul{W}^-(s)=r\},\\
\gb_r^-&:=-\min\{s\ge0:W^-(s)=\ul{W}^-(\tau_r^-)\}.
\end{align*}
And similarly define $\tau_r^-, \gb_r^-$ for $W^-$. One can see
that with probability one, it holds $b_r\in\{\gb_r^-, \gb_r^+\}$
(see e.g. \cite{ZE}). In the following, we will use the fact that
$-W(\gb_1^+)$ is exponential random variable with mean 1 (see
lemma of \S 1 in \cite{NP}).

Also we will use the following fact (immediate consequence of
Lemma 1.1.1 in \cite{CR}).

\medskip

\nid \textbf{Fact 1}: There is a constant $C$ so that for any
$\rho>0, h\in[0,\rho], v>0$, it holds
$$ \PPP(\sup_{y,z\in[0, \rho], |y-z|\le h} |W(y)-W(z)|\ge
v\sqrt{h})\le C \frac{\rho}{h}e^{-v^2/3}.
$$

\begin{lemma} \label{LILbound}
For all $c>1, M>1$, and $a>0$, it holds
$$\PPP(\sup_{1\le t\le M} \frac{|B_t|}{\sqrt{t}}\ge a)\le
4 (\frac{\log M}{\log c}+1) \frac{c}{a}
e^{-\frac{1}{2}(\frac{a}{c})^2}.
$$
\end{lemma}

\begin{proof}
Let $N:=[\frac{\log M}{\log c}]$, and $t_n:=c^n$ for $n=0,\ldots
N+1$. Then
\begin{multline*}
\PPP(\sup_{1\le t\le M} \frac{|B_t|}{\sqrt{t}}\ge a)\le
\sum_{n=0}^N \PPP(\sup_{t_n\le s\le t_{n+1}}
\frac{|B_s|}{\sqrt{s}}\ge a) \le \sum_{n=0}^N \PPP(\sup_{0\le s\le
t_{n+1}}|B_s|\ge \sqrt{t_n} a)\\ \le 4 \sum_{n=0}^N \PPP(B_1\ge
\frac{a}{c})\le 4 (\frac{\log M}{\log c}+1) \frac{c}{a}
e^{-\frac{1}{2}(\frac{a}{c})^2}
\end{multline*}
\end{proof}

\begin{lemma} \label{tailslemma} For all $x>0$, we have

1. $\PPP(W^{\#}(b_{1^-}, b_1)<x)\le 4 \sqrt{x}.$

\medskip

2. $\PPP(W^{\#}(b_1, b_{1^+})>x)< 2(2+x) e^{-x}.$

\medskip

3. $\PPP(W^{\#}(b_1, 0)<x)\le 2x.$

\medskip

4. $\PPP(W^{\#}(b_1, 0)>x)< 6 e^{-x}.$

\medskip

5. $\PPP(|b_1|>x)< 2e^{-x}$.

6. $(W(b_1)-W(b_{1^++}))/W^{\#}(b_1, b_{1^++})$ is an exponential
random variable with mean 1.

\end{lemma}
\begin{proof}

1. $W^{\#}(b_{1^-}, b_1)\ge \min\{\ol{W}(\gb^+_1),
\ol{W}(\gb^-_1)\}$. So that $\PPP(W^{\#}(b_{1^-}, b_1)<x)\le 2
\PPP(\ol{W}(\gb^+_1)<x)<4\sqrt{x}$ by Lemma \ref{exrursionlemma}.

2. $W^{\#}(b_1, b_{1^+})=\min\{h_1, h_2\}$, the min of the two
$1$-slopes with bottom $b_1$. The one has density
$\textbf{1}_{x\ge 1}e^{-x+1}$, and the other, $\textbf{1}_{x\ge 1}
(2x-1)e^{-x+1}/3$. So that $\PPP(W^{\#}(b_1, b_{1^+})>x)\le
\PPP(h_1>x)+\PPP(h_2>x)$. An easy computation shows that, for
$x>1$, the last quantity equals $e^{-x+1}+(1+2x) e^{-x+1}/3<(4+2x)
e^{-x}.$

3. $W^{\#}(b_1,0)\ge \min\{-W(\gb_1^+), -W(\gb_1^-)\}$. So that
$\PPP(W^{\#}(b_1, 0)<x)\le 2\PPP(-W(\gb_1^+)<x)<2x$, since
$-W(\gb_1^+)$ has exponential distribution with mean 1.

4. $W^{\#}(b_1, 0)=\ol{W}(b_1)-W(b_1)<1-W(b_1)$, and $-W(b_1)\le
\max\{-W(\gb_1^+), -W(\gb_1^-)\}$. So that $\PPP(W^{\#}(b_1,
0)>x)\le \PPP(-W(b_1)
>x-1)\le 2 \PPP(-W(\gb_1)
>x-1)\le 2 e^{1-x}$.

5. The density of $b_1$ is
$f_{b_1}(x):=\frac{2}{\pi}\sum_{k=0}^{+\infty} \frac{(-1)^k}{2k+1}
e^{-(2k+1)^2|x|}$, and the required inequality follows after
integration.

6. For $\ell\in\D{R}$, let
\begin{align*}
H_\ell^-:=& \sup \{s<0: W_s=\ell\},\\
H_\ell^+:=& \inf \{s>0: W_s=\ell\}, \\
K_\ell:=&\min \{\max_{H_\ell^-\le s \le 0}W_s, \max_{0\le s \le
H_\ell^+}W_s\}.
\end{align*}
Also let $\ell_0:=\sup\{\ell<0: \text{ one of $H_\ell^-, H_\ell^-$
jumps at $\ell$, and } K_\ell+\ge 1\}$. Clearly, $\ell_0=W(b_1)$.
Assume that $\max_{H_{\ell_0}^-\le s \le 0}W_s>\max_{0\le s \le
H_{\ell_0}^+}W_s$. Then $b_{1^++}>0$, and $W^{\#}(b_1,
b_{1^++})=\max_{0\le s \le H_\ell^+}W_s$.  The way to locate
$(b_{1^++}, W(b_{1^++})$ is as follows. Let $H_{\ell_0+}^+$ denote
the left limit of $H_{\ell}^+$ at $\ell_0$. We condition on
$W|[H_{\ell_0}^-, H_{\ell_0+}^+]$, we look at
$B=\{W(H_{\ell_0+}^++s)-W(b_r): s>0\}$ (which is a standard
Brownian motion), and we wait until $B-\ul{B}$ hits $W^{\#}(b_1,
b_{1^++})$. When this happens, the value of $-\ul{B}$ is an
exponential random variable with mean $W^{\#}(b_1, b_{1^++})$ (see
beginning of this section). This proves our claim.
\end{proof}

\begin{lemma} \label{denslemma}
 $x_1(W,1), |x_0(W,1)|$ have densities bounded by 1.
\end{lemma}

\begin{proof}
It is true that $\{x_{k+1}(W,1)-x_k(W,1):
k\in\D{Z}\setminus\{0\}\}$ is a set of i.i.d. random variables
with the same distribution as $\ell:=\inf\{s\ge0:|W_s|=1\}$ (see
proposition of \S 1 in \cite{NP}) . Call $f_\ell(x)$ the density
of this random variable. Since for any fixed $t$, the process
$(W_{s-t}-W_{-t}:s\in\D{R})$ is a standard Brownian motion, one
can take $t\to+\infty$ and use the renewal theorem to show that
$x_1(W,1), |x_0(W,1)|$ are respectively the residual waiting time
after 0 and the age at time 0 for a renewal process ``starting at
$-\infty$'' and with increments having distribution $\ell$. Their
densities are computed in Exercise 4.7 of Chapter 3 in \cite{DU},
and they both equal $\int_x^{+\infty} f_\ell(z)dz$.
\end{proof}

\nid \textbf{Fact 2}: Let $W$ be standard Brownian motion. For the
time $\rho(1):=\inf\{s>0:W(s)=1\}$, it holds
$$\PPP(\rho(1)>u)<u^{-1/2} \text{ for all } u>0.$$

\nid This follows from
$\PPP(\rho(1)>u)=\PPP(\ol{W}(u)<1)=\PPP(\ol{W}(1)<u^{-1/2})$ and
the fact that $\ol{W}(1)$ has density $\sqrt{2/\pi}
e^{-x^2/2}1_{x\ge0}$.


\medskip

\begin{lemma} \label{thm2lemma}

For the three choices of $y_0, \ga, \gb, \gga, \gz, \eta, H$ in
Lemma \ref{locprop} and the two in Lemma \ref{transition}, we have
for large $r$

\medskip

1. $\PPP(|\eta-y_0|> H^4)<r^{-1/4}$.

\medskip

2. $\PPP(\max_{s\in[y_0, \eta]}W(s)>H^{3/2})<cr^{-1/4}.$

\medskip

3. $\PPP(A(\eta)/(A(\eta)+|A(y_0)|)>H^{-2})<r^{-2}.$

\medskip

4. $\PPP(\frac{A(\eta)-A(\gb)}{A(\gz)-A(\gb)}>e^{H^2})<r^{-1/3}.$

\medskip

5. $\PPP(A(\gb)/A(\gz)<-H^{-5})<Cr^{-1/2}$.

\medskip

 6. $\PPP(\sup_{y\in[\gz, \eta]} \sup_{y_0<s<y, s\notin(\ga, \gga)}
\int_s^{y} e^{W(t)-W(s)} dt/\int_{\gb}^{y}e^{W(t)-W(\gb)} dt >
H^{-\min\{ k_1, k_3\}+5})<c r^{-1/4}$.

\medskip

7. $\PPP(\sup_{y\in[\gz, \eta]} \sup_{y_0<s<y, s\notin(\ga, \gga)}
(e^{W(s)}\int_0^y e^{W(t)-W(\gb)}dt)^{-1}
>\exp(H^{3/4}))<7r^{-1/4}.$
\end{lemma}

\begin{proof}

In all uses, it is $H=r$.

\medskip

\nid 1. We prove this claim at once for all the cases that we use
it. Let
\begin{align*}\rho(r)&:=\inf\{s>0:W(s)=r\},\\
\tau_1(r)&:=\inf\{s>\rho(r): W(s)=\ul{W}(\rho(r))\},\\
\tau_2(r)&:=\inf\{s>\tau_1(r): W(s)=\ol{W}(\tau_1(r))\}.
\end{align*}
 Through
the analogous series of definitions we define $\tilde \tau_2(r)$
for the path $(W(-\cdot))_{s\ge0}$. In all cases that we use the
lemma, it holds $[y_0, \eta]\subset [-\tilde\tau_2(r),
\tau_2(r)]$. Also let $\tau_0(r):=\inf\{s\in[0,
\tau_1(r)]:W(s)=\ol{W}(\tau_1(r))\}.$

Clearly, $\tau_2(r)\overset{law}{=}r^2\tau_2(1)$. We write
$$
\tau_2(1)=\rho(1)+\tau_0(1)-\rho(1)+\tau_2(1)-\tau_0(1).$$
First, $\PPP(\rho(1)>r)<r^{-1/2}$ for all $r>0$.

 The
random variable $r_1:=W(\tau_0(1))-1/(1-\ul{W}(\rho(1)))$ has
density $(1+x)^{-2}1_{x\ge0}$ because
\begin{multline*}
\PPP(r_1>x)= \EEE\{\PPP((W(\tau_0(1))-1>x(1-\ul{W}(\rho(1)))
|\ul{W}(\rho(1))))\}
\\= \EEE\{\PPP(W \text{ starting from $1$ hits first
$1+x(1-\ul{W}(\rho(1))$ and then } \ul{W}(\rho(1))
|\ul{W}(\rho(1))))\}\\
=\EEE(\frac{1-\ul{W}(\rho(1))}{1-\ul{W}(\rho(1))+x(1-\ul{W}(\rho(1)))})=(1+x)^{-1}.
\end{multline*}
We used the Markov property on the stopping time $\rho(1)$. Now
given the values of $\ul{W}(\rho(1)), W(\tau_0(1))-1$, the law of
$\tau_0(1)-\rho(1)$ is the same as the time it takes for a three
dimensional Bessel starting from $1+\ul{W}(\rho(1))$ to hit
$W(\tau_0(1))+\ul{W}(\rho(1))$. So it is bounded stochastically
from above by the time it takes for Brownian motion starting from
zero to hit $W(\tau_0(1))-1$. This last time equals in law to
$(W(\tau_0(1))-1)^2 X=r_1^2 (1-\ul{W}(\rho(1)))^2 X$, where $r_1$
has density $(1+x)^{-2}1_{x\ge0}$, $X$ has the same law as
$\rho(1)$, and $X, r_1, \rho(1)$ are independent. Also it is easy
to see that $1-\ul{W}(\rho(1)$ has the same distribution as $r_1$.

\nid Consequently
$$
\PPP(\tau_0(1)-\rho(1)>r^{3/2})\le
\PPP(r_1>r^{1/4})+\PPP(1-\ul{W}(\rho(1))>r^{1/4})+\PPP(X>r^{1/2})<3r^{-1/4}.
$$
As above, we show that the random variable $r_2:=
(W(\tau_1(1))-\ul{W}(\tau_2(1)))/(W(\tau_0(1))-\ul{W}(\rho(1)))$
has density $(1+x)^{-2}1_{x\ge0}$. Given
$W(\tau_0(1))-\ul{W}(\tau_2(1))$, the law of $\tau_2(1)-\tau_0(1)$
is the same as the lifetime of a Brownian excursion conditioned to
have height $W(\tau_0(1))-\ul{W}(\tau_2(1))$. This is equal in law
to $(W(\tau_0(1))-\ul{W}(\tau_2(1)))^2(Y_1+Y_2)$, where $Y_1, Y_2$
have the law of the time it takes for a three dimensional Bessel
process starting from zero to hit 1. Observe that
\begin{multline*}
W(\tau_0(1))-\ul{W}(\tau_2(1))=W(\tau_0(1))-W(\tau_1(1))+W(\tau_1(1))-\ul{W}(\tau_2(1))\\=
W(\tau_0(1))-W(\tau_1(1))+r_2(W(\tau_0(1))-W(\tau_1(1))) \\ =
(1+r_2)(W(\tau_0(1))-W(\tau_1(1)))=(1+r_2)(1+r_1)(1-\ul{W}(\rho(1)))
\end{multline*}
So that
$$\tau_2(1)-\tau_0(1)=
(1+r_2)^2(1+r_1)^2(1-\ul{W}(\rho(1)))^2(Y_1+Y_2),$$
and
\begin{multline*}
\PPP(\tau_2(1)-\tau_0(1)>r^2/2)\le
\PPP(1+r_1>r^{1/4})+\PPP(1+r_2>r^{1/4})\\ +
\PPP(1-\ul{W}(\rho(1))>r^{1/4})
+\PPP(Y_1>r^{1/2})+\PPP(Y_2>r^{1/2})<5r^{-1/4}.
\end{multline*}
Of course, $\PPP(\rho(1)>r^{1/2})<r^{-1/4}$. Combining all the
above estimates, we get that on a set whose complement has
probability at most $9r^{-1/4}$, it holds $\tau_2(1)\le
r^{1/2}+r^{3/2}+r^2/2$, which is less than $r^2$ for $r>9$. This
finishes the proof of part 1.

\bigskip

\nid 2. As we mentioned in the proof of part one, in all uses of
the lemma, it holds $[y_0, \eta]\subset [-\tilde\tau_2(r),
\tau_2(r)]$. So $\max_{s\in[y_0, \eta]}W(s)\le \max\{W(\tau_2(r)),
W(-\tilde \tau_2(r))\}=\max\{W(\tau_0(r)), W(-\tilde
\tau_0(r))\}$. Now $W(\tau_0(r))\overset{law}{=}r W(\tau_0(1))$,
and we saw that
$W(\tau_0(1))=1+W(\tau_0(1))-1=1+r_1(1-\ul{W}(\rho(1)))$. Since
$$\PPP(1+r_1>r^{1/4})=\PPP(1-\ul{W}(\rho(1))>r^{1/4})=r^{-1/4},$$
for $r>1$ we have outside a set of probability at most
$2r^{-1/4}$, that $W(\tau_0(1))<1+(r^{1/4}-1)r^{1/4}<r^{1/2}$.
Consequently, $\PPP(\max\{W(\tau_2(r)), W(\tilde
\tau_2(r))\}>r^{3/2})\le 4r^{-1/4}$.

\bigskip

\nid 7. On $[\gz_r<1]^c$, we have
$$\frac{e^{-W(s)}}{\int_0^y e^{W(t)-W(b_r)}dt}\le
\exp(-\inf_{x_0<s<\eta_r}W(s)).$$
Recall the definitions made above, in the proof of part 1. The
exponent of the last expression is bounded above by
 $\max\{-\ul{W}(\tau_2(r)),
-\ul{W}(\tau_2(r))\}$, which has the same distribution as
$r\max\{-\ul{W}(\tau_2(1)), -\ul{W}(\tau_2(1))\}$. Observe that
\begin{multline*}
-\ul{W}(\tau_2(1))=-\ul{W}(\rho(1))+
W(\tau_1(1))-\ul{W}(\tau_2(1))
=-\ul{W}(\rho(1))+r_2(W(\tau_0(1))-\ul{W}(\rho(1)))\\ =
-\ul{W}(\rho(1))+r_2(W(\tau_0(1))-1+1-\ul{W}(\rho(1))) =
-\ul{W}(\rho(1))+r_2(1+r_1)(1-\ul{W}(\rho(1))),
\end{multline*}
and $\PPP(\max\{r_2, 1+r_1,
1-\ul{W}(\rho(1))\}>r^{1/4})<3r^{-1/4}$. On the complement of
$[\max\{r_2, 1+r_1, 1-\ul{W}(\rho(1))\}>r^{1/4}]$, we have
$-\ul{W}(\tau_2(1))<r^{1/4}+r^{3/4}$. Consequently, for $r>1$,
$$\PPP(-\inf_{x_0<s<\eta_r}W(s)>r^{7/4})<6r^{-1/4}.$$
 Also $\PPP(\gz_r<1)<r^{-2}$, and this proves the statement.

\medskip

\nid The events
\begin{align*}
A_4:&=[x_0(W,r^++)\in(-1,0]]\cup [|x_0(W,r^++)|>r^4]\cup
[\gz_r<1]\cup [|\gz_r|>r^4]\cup [|\eta_r|>r^4], \\
A_5:&= [\sup_{y,z\in[-r^4, r^4], |y-z|\le 1} |W(y)-W(z)|\ge \log
r],
\end{align*}
will be used bellow. Observe that
$$\PPP(x_0(W,r^++)\in(-1,0])\le
\PPP(x_0(W,r)\in(-1,0])=\PPP(x_0(W,1)\in(-1/r^2,0])\le r^{-2},$$
and $\PPP(\gz_r<1)=\PPP(\gz_1<r^{-2})\le
\PPP(x_1(W,1)<r^{-2})<r^{-2}$. The last inequality follows from
Lemma \ref{denslemma}. Combining these with part 1 of the lemma,
we get $\PPP(A_4)<cr^{-2}$. Also, applying Fact 1, with $h=1,
\rho=r^4, v=\log r/2$, we get
$$\PPP(A_5)\le C
r^4e^{-(\log r)^2/12}<c r^{-2}.$$

\nid 5. In four of the five cases we use the lemma, it holds
$A(b_r)/A(\gz_r)>0$, and we have nothing to prove. The only case
where something needs a proof is in the claim
$P_W(F_{X^{W^*}}([\tau(-\tilde \gz_r), \tau(-\eta_r)])\nsubseteq
(-\gga_r, -\ga_r))<cr^{-1/4}$ contained in the proof of Lemma
\ref{locprop} (Case 2). Let $A_7:=[W^{\#}(b_r,0)<\sqrt{r}]\cup
[|b_r|>r^3]\cup[\tilde \gz_r>-1]$. Then $\PPP(A_7)\le
\PPP(W^{\#}(b_1,0)<1/\sqrt{r})+\PPP(|b_1|>r)+\PPP(\tilde \gz_1
>-r^{-2})\le 2/\sqrt{r}+2e^{-r}+cr{-1/2}<Cr^{-1/2}$. On $(A_4\cup A_5 \cup A_7)^c$ we
have
$$\frac{A(b_r)}{|A(\tilde \gz_r)|}=\frac{\int_0^{b_r}
e^{W(y)-W(b_r)}dy}{\int_{\tilde \gz_r}^0 e^{W(y)-W(b_r)}dy}\le
\frac{b_r e^{W^{\#}(b_r,0)}}{e^{W^{\#}(b_r,0)+k_2\log
W^{\#}(b_r,0)-\log r}}=\frac{rb_r}{(W^{\#}(b_r,0))^{k_2}}\le
r^{4-k_2/2}<r^{-5}$$
since $k_2\ge18$.

\bigskip

The remaining parts of the lemma we prove them only for the first
choice of $y_0, \ga, \gb, \gga, \gz, \eta, H$, i.e., $(y_0, \ga,
\gb, \gga, \gz, \eta, H):=(x_0(W,r^++), \ga_r, b_r, \gga_r, \gz_r,
\eta_r, r)$. For the others, the proof is similar.

\medskip

\nid 3. The quotient inside the probability equals
$$\frac{A(\eta_r)}{A(\eta_r)+
|A(x_0(W,r^+))|}=\frac{1}{1+|A(x_0(W,r^+))|/A(\eta_r)},
$$
and we will show that $|A(x_0(W,r^+))|/A(\eta_r)$ is large. We
will use $x_0$ instead of $x_0(W,r^+)$ in the following. On
$(A_4\cup A_5)^c$ we have $\sup_{s\in[x_0, x_0-1]}
|W(s)-W(x_0)|<\log r$. So that $W(s)\ge W(x_0)-\log r$ on $[x_0,
x_0-1]$, and
\begin{multline*}
\frac{|A(x_0(W,r^+))|}{A(\eta_r)}=\frac{\int_{x_0}^0e^{W(y)}dy}{\int_0^{\eta_r}e^{W(y)}dy}
=
\frac{\int_{x_0}^0e^{W(y)-W(b_r)}dy}{\int_0^{\eta_r}e^{W(y)-W(b_r)}dy}\ge
\frac{e^{W(x_0)-W(b_r)-\log r}}{\eta_r e^{W^{\#}(b_r,b_{r^+})-k_3
\log W^{\#}(b_r,b_{r^+})}} \\ \ge \frac{e^{k_3 \log
W^{\#}(b_r,b_{r^+})-\log r}}{r^4}=
\frac{(W^{\#}(b_r,b_{r^+}))^{k_3}}{r^5}\ge  r^{k_3-5}.
\end{multline*}
We used that fact that $W(x_0)-W(b_r)\ge W^{\#}(b_r,b_{r^+})$,
which holds because we assumed that $b_{r^++}>0$, and also that
$W^{\#}(b_r,b_{r^+})\ge r$.
 Thus, on the complement of $A_4\cup A_5$, it holds
$$\frac{A(\eta_r)}{A(\eta_r)+ |A(x_0(W,r^+))|}\le \frac{1}{1+r^{k_3-5}}\le  r^{-2}.$$
since $k_3\ge 7$.
\medskip

\nid 4. Let $A_6:=[W^{\#}(b_{r^-}, b_r)<\log r]\cup[W^{\#}(b_r,
b_{r^+})>r^{3/2}]$. Then $$\PPP(A_6)=\PPP(W^{\#}(b_{1^-},
b_1)<\log r/r)+ \PPP(W^{\#}(b_1, b_{1^+})>\sqrt{r}) < r^{-1/3},$$
using Lemma \ref{tailslemma}. Now on $(A_4\cup A_5\cup A_6)^c$ we
have $\gz_r-b_r>1$ (because of the definition of $A_5$ and the
fact that $W(\gz_r)-W(b_r)>W^{\#}(b_{r^-}, b_r)>\log r$ on
$A_6^c$), and
\begin{multline*}
\frac{A(\eta_r)-A(b_r)}{A(\gz_r)-A(b_r)}=\frac{\int_{b_r}^{\eta_r}
e^{W(y)-W(b_r)}dy}{\int_{b_r}^{\gz_r} e^{W(y)-W(b_r)}dy}\le
\frac{\eta_r e^{W^{\#}(b_r, b_{r^+})-k_3\log W^{\#}(b_r,
b_{r^+})}}{e^{W^{\#}(b_{r^-}, b_r)+k_2\log W^{\#}(b_{r^-}, b_r)
-\log r}}=\frac{r\eta_r e^{ W^{\#}(b_r, b_{r^+})-W^{\#}(b_{r^-},
b_r)}}{W^{\#}(b_r, b_{r^+})^{k_3} W^{\#}(b_{r^-}, b_r)^{k_2}} \\
\le \frac{r^{5-k_3} e^{ W^{\#}(b_r, b_{r^+})-W^{\#}(b_{r^-},
b_r)}}{W^{\#}(b_{r^-}, b_r)^{k_2}}<r^{5-k_3} e^{r^{3/2}}<e^{r^2}.
\end{multline*}
The last inequality holding for big $r$. Finally note that
$\PPP(A_4\cup A_5\cup A_6)<r^{-1/3}.$

\medskip

\nid 6. The quantity of interest is
$$\frac{\int_s^y
e^{W(t)-W(s)} dt}{\int_{b_r}^y e^{W(t)-W(b_r)} dt}.$$
Let $m_r:=\inf \{x>b_r:W(x)-W(b_r)=W^\# (b_r, b_{r^++})\}$ and
$$A_7:=[\text{ there is an } s\in[x_0(W,r^++), m_r]\setminus (\ga_r,
\gga_r)\text{ with } W(s)-W(b_r)\le k_1\log r]$$
 For
$s\in[j_r,l_r]\setminus (\ga_r, \gga_r)$ it holds $W(s)-W(b_r)\ge
k_1\log r$ by the definition of $\ga_r, \gga_r$. It remains to
study the intervals $[x_0(W,r^++), j_r]$, $[l_r, m_r]$.
 We will study only the first, the case of the second is similar.
$B:=(W(-s+j_r)-W(j_r))_{s\ge0}$ is a standard Brownian motion. If
there is $s\in[x_0(W,r^++), j_r]$ with $W(s)-W(b_r)<k_1\log r$,
then $B$ visits $-r+k_1\log r$ and then returns to 0 before
hitting $-r$. This last event has probability $k_1\log r/r$. So
that $\PPP(A_7)<r^{-1/2}$ for large $r$. Finally,
$\PPP(m_r>r^4)<r^{-1/4}$ from part 1 of the lemma.

Working as in part 4, we see that on $([m_r>r^4]\cup[\gz_r<1]\cup
A_5 \cup A_7)^c$ we have the following bounds.

\nid If $y<m_r$, then the bound is
$$\frac{(y-s)\exp(W^{\#}(b_r, y)-k_1\log r)}{\exp(W^{\#}(b_r,
y)-\log r)}<(\eta_r-x_0(W,r)) r^{-k_1+1}<r^{-k_1+5}.$$
If $y>m_r$, then the bound is
$$\frac{(\eta_r-x_0(W,r)) \exp(W^{\#}(b_r, b_{r^+})-k_3\log W^{\#}(b_r, b_{r^+}))}{\exp(W^{\#}(b_r, b_{r_1+})-\log r)}<
(\eta_r-x_0(W,r)) (W^{\#}(b_r, b_{r^+}))^{-k_3}r <r^{-k_3+5}.$$
We used the definition of $\eta_r$ to bound the numerator.

\end{proof}

\begin{lemma} \label{hitlemma}
For the sets $A_3, \gS_3$ defined in lemmata \ref{locprop},
\ref{transition} respectively, it holds
$$\C{P}(A_3^c)<C r^{-1/4}, \C{P}(\gS_3^c)<C r^{-1/4}.$$
\end{lemma}

\begin{proof}
We have
\begin{multline*}
\C{P}(A_3^c)\le \C{P}(\tau(\tilde
\gz_r)<\tau(\gz_r))+\C{P}(\tau(\hat
\gz_r)<\tau(\tilde \gz_r))+\C{P}(\tau(\eta_r)<\tau(\hat \gz_r))\\
= \EEE(\frac{A(\gz_r)}{A(\gz_r)+|A(\tilde
\gz_r)|})+\EEE(\frac{|A(\tilde \gz_r)|}{|A(\tilde
\gz_r)|+A(\hat\gz_r)})+\EEE(\frac{A(\hat\gz_r)}{A(\hat\gz_r)+|A(\eta_r)|})\end{multline*}
We work as in part 3 of Lemma \ref{thm2lemma}. Let $A_5$ be
defined as there, and
\medskip

$$A_8:=[\gz_r<1
\text{ or } \tilde \gz_r\in(-1, 0] \text{ or } \hat \gz_r>r^4
\text{ or } |\eta_r|
>r^4 \text{ or } |\gz_r|>r^4 \text{ or } W^{\#}(b_r,0)>\sqrt{r} \text{ or }
|\ol{W}(\gb_r^+)-\ol{W}(\gb_r^-)|<\sqrt{r}]. $$

\medskip

\nid Then $\PPP(\gz_r<1 \text{ or } \tilde \gz_r\in(-1, 0])\le
\PPP(x_1(W,r)<1)+\PPP(\tilde \gz_r\in(-1, 0])<r^{-2}+cr^{-1/2}$
using lemmata \ref{denslemma}, \ref{ztildlemma}. So that
$\PPP(A_8)\le r^{-2}+cr^{-1/2}+ c'
r^{-1/4}+2r^{-1/2}+3r^{-1/4}<Cr^{-1/4}$.

\nid On $(A_8\cup A_5)^c$ the quantities $|A(\tilde
\gz_r)|/A(\gz_r)$, $A(\hat\gz_r)/|A(\tilde \gz_r)|$,$
|A(\eta_r)|/A(\hat\gz_r)$ are large. Indeed
\begin{align*}\frac{|A(\tilde \gz_r)|}{A(\gz_r)}&=\frac{\int_{\tilde
\gz_r}^0 e^{W(y)-W(b_r)}dy}{\int_0^{\gz_r} e^{W(y)-W(b_r)}dy} \ge
\frac{\exp(W^{\#}(b_r,0)+2 k_2 \log W^{\#}(b_r,0)-\log r)}{\gz_r
\exp(W^{\#}(b_{r^-}, b_r)+k_2 \log W^{\#}(b_{r^-},b_r))}\ge
\frac{(W^{\#}(b_r,0))^{k_2}}{r^5}\ge
r^{k_2/2-5},\\
\frac{A(\hat \gz_r)}{|A(\tilde \gz_r)|}&=\frac{\int_0^{\hat \gz_r}
e^{W(y)-W(b_r)}dy}{\int_{\tilde \gz_r}^0 e^{W(y)-W(b_r)}dy}  \ge
\frac{\exp(W^{\#}(b_r,0)+3k_2 \log W^{\#}(b_r,0)-\log r)}{|\tilde
\gz_r| \exp(W^{\#}(b_r,0)+2k_2\log W^{\#}(b_r,0))}\ge
\frac{(W^{\#}(b_r,0))^{k_2}}{r^5}\ge r^{k_2/2-5}, \\
\frac{|A(\eta_r)|}{A(\hat \gz_r)}&=\frac{\int_{\eta_r}^0
e^{W(y)-W(b_r)}dy}{\int_0^{\hat \gz_r} e^{W(y)-W(b_r)}dy}  \ge
\frac{\exp(W^{\#}(b_r,b_{r^+})-\log r)}{\hat \gz_r
\exp(W^{\#}(b_r,0)+3k_2\log W^{\#}(b_r,0))}=
\frac{\exp(W^{\#}(b_r,b_{r^+})-W^{\#}(b_r, 0))}{\hat \gz_r r
(W^{\#}(b_r, 0))^{3k_2}}.
\end{align*}
In the first line, we used the fact that $W^{\#}(b_{r^-}, b_r)\le
W^{\#}(b_r,0)$. Regarding the last quantity of the third line,
observe that $W^{\#}(b_r,b_{r^+})-W^{\#}(b_r, 0)\ge
|\ol{W}(\gb_r^+)-\ol{W}(\gb_r^-)|\ge \sqrt{r}$, and $\hat \gz_r r
(W^{\#}(b_r, 0))^{3k_2}\le r^{5+3k_2/2}$. So that
$|A(\eta_r)|/A(\hat \gz_r)\ge e^{\sqrt{r}}r^{-5-3k_2/2}>r^2$ for
large $r$.

\nid Finally,
$$\C{P}(\gS_3^c)=\EEE(\frac{A(\gz_{r^+})-A(\eta_r)}{A(\gz_{r^+})-A(z)}).$$
The quantity in the expectation is always at most one. The set
$A_5\cup [\gz_{r^++}<1 \text{ or } \gz_{r^++}>r^4]$ has
probability at most $9r^{-1/4}+r^{-2}$ (because of
$0<x_1(W,r)<\gz_{r^++}$ and Lemma \ref{denslemma}), and on its
complement, it holds
\begin{multline*}\frac{A(\gz_{r^+})-A(\eta_r)}{A(\gz_{r^+})-A(z)}=
\frac{\int_{\eta_r}^{\gz_{r^+}}
e^{W(y)-W(b_{r^+})}dy}{\int_z^{b_{r^+}} e^{W(y)-W(b_{r^+})}dy} \le
\frac{\gz_{r^+}
\exp(W(\gz_{r^+})-W(b_{r^+}))}{\exp(W(\gz_{r^+})-W(b_{r^+})+k_2\log
r-\log r)}\\ =\gz_{r^+}/r^{k_2-1}<r^{5-k_2}<r^{-2}.
\end{multline*}

\end{proof}

\begin{lemma} \label{exrursionlemma}
The random variable $\ol{W}(\gb_1^+)$ has density
$f(x)=-\textbf{1}_{x\in(0,1]}\log x$. In particular,
$\PPP(\ol{W}(\gb_1^+)<x)<2\sqrt{x}$ and
$\PPP(|\ol{W}(\gb_1^+)-\ol{W}(\gb_1^-)| <x)<3\sqrt{x}$ for all
$x\in[0,1]$.
\end{lemma}
\begin{proof}
The proof uses excursion theory, for which we give the basic
setup. The following are standard.

Consider the process $Y(t):=W(t)-\ul{W}(t)$. A local time process
for $Y$ is $-\ul{W}$. Let $(\gep_t)_{t>0}$ be corresponding
excursion process. For any $\gep$ in the space of the excursions,
we denote by $\ol{\gep}$ the maximum value of $\gep$. The process
$\{(\ol{\gep}_t, t):t\ge0\}$ is a Poisson point process in
$[0,+\infty)\times [0,+\infty)$ with characteristic measure
$x^{-2} dx dt$. The time $t^*:=\inf\{s\ge0: \ol{\gep}_s\ge1\}$ has
exponential distribution with mean one and $\{(\ol{\gep}_t, t):
t\le t^*\}$ has the same law as the restriction in $[0,1]\times
[0,\tau]$ of a Poisson point process in $[0,1]\times [0,+\infty)$
with characteristic measure $dn:=x^{-2}dxdt$, where $\tau$ is an
exponential random variable independent of the process. Let $N$ be
the counting measure of that process. Also for all $t>0$, let
$A_t:=\{(y,s):s\in[0,t] \text{ and }  y>x+s\}$. Then
$n(A_t)=\int_0^{\min\{t, 1-x\}} \int_{s+x}^1 y^{-2} dy ds$, which
equals $\log(1+t/x)-t$ if $t<1-x$, and $-\log x-1+x$ otherwise.
Then
\begin{multline*}
\PPP(\ol{W}(\gb_1^+)<x)=\PPP(\text{ for all } s<t^*, \text{ it
holds } \ol{\gep}_s-s<x)= \int_0^{+\infty} e^{-t} P(N(A_t)=0) dt=
\int_0^{+\infty} e^{-t} e^{-n(A_t)}dt\\=\int_0^{1-x}
\frac{x}{x+t}dt+\int_{1-x}^{+\infty} e^{-t}xe^{1-x} dt=-x\log x+x.
\end{multline*}
In particular, the density is $f(x)=-\log x$.  To bound
$\PPP(\ol{W}(\gb_1^+)<x)$, we observe that for $x\in(0,1]$, it
holds $-x\log x+x<2\sqrt{x}.$ To bound
$\PPP(|\ol{W}(\gb_1^+)-\ol{W}(\gb_1^-)| <x)$, we use the fact that
$f$ is decreasing in $(0,1]$ to get
$$
\PPP(|\ol{W}(\gb_1^+)-\ol{W}(\gb_1^-)| <x)=\int_0^1 f(y)
\PPP(\ol{W}(\gb_1^+)\in(y-x,y+x))dy \le \int_0^1 f(y)
\PPP(\ol{W}(\gb_1^+)<2x)dy<2\sqrt{2x}.
$$
\end{proof}

\begin{lemma} \label{ztildlemma}
For all $x>0$, it holds
$$\PPP(\tilde \gz_1\in[-x,0])\le c x^{1/4}.$$
\end{lemma}

\begin{proof}
Remember the definitions in the beginning of this section, and let
$\rho_{W^-}(c)=\inf\{t>0: W^-(t)=c\}$.
\begin{multline*}
\PPP(\tilde \gz_1\in[-x,0])\le
\PPP(\rho_{W^-}(\ol{W}(\gb_1^+))<x)=-\int_0^1
\PPP(\rho_{W^-}(s)<x)
\log s ds=-\int_0^1 \PPP(\rho_{W^-}(1)<x/s^2) \log s ds\\
=\int_1^{+\infty} \frac{\log y}{y^2}\PPP(\rho_{W^-}(1)<xy^2) dy<
\int_1^{x^{-1/2}} \frac{\log y}{y^2}(xy^2)^{1/2}
dy+\int_{x^{-1/2}}^{+\infty} \frac{\log y}{y^2} dy\\
=\sqrt{x}\int_1^{x^{-1/2}} \frac{\log
y}{y}dy+\sqrt{x}(1-\frac{\log x}{2})
\end{multline*}
We used Lemma \ref{exrursionlemma} for the density of
$\ol{W}(\gb_1^+)$, and Fact 2. The last quantity is easily shown
to have bound of the form $cx^{1/4}$.
\end{proof}


The next lemma says that, with high probability, the points
$\tilde z_r, \hat \gz_r, m_r, z$ are are as we depict them in
Figures \ref{Fig3}, \ref{Fig4}. Parts 1 and 2 should be used when
one proves the versions of Lemma \ref{thm2lemma} needed in the
proof of Lemma \ref{locprop}. Part 3 is used in the proof of Lemma
\ref{transition}.

\begin{lemma} \label{locationslemma}

1. $\PPP(b_{r^++}>\tilde \gz_r)<cr^{-1/2}$.

2. $\PPP(\min_{s\in[\gz_r, \hat \gz_r]}W(s)>W(b_r))<cr^{-1/2}$

3. $\PPP(z<m_r)<c r^{-1/2}$.
\end{lemma}

\begin{proof}

Remember the definitions in the beginning of this section.

 1. The probability of interest is
bounded by twice the following probability (since $b_r$ will be
either $\gb_r^+$ or $\gb_r^-$.)
\begin{multline*}\PPP(W^-\text{ starting from $\ol{W}(\gb_r^+)$ hits $W(\gb_r^+)$
before }\ol{W}(\gb_r^+)+ k_2 \log
(\ol{W}(\gb_r^+)-W(\gb_r^+)))\\=\EEE(\frac{2k_2 \log
(\ol{W}(\gb_r^+)-W(\gb_r^+))}{\ol{W}(\gb_r^+)-W(\gb_r^+)}<4k_2\EEE(\frac{1}{\sqrt{
\ol{W}(\gb_r^+)-W(\gb_r^+)}})< \frac{4k_2}{\sqrt{r}}
\EEE(\frac{1}{\sqrt{\ol{W}(\gb_1^+)}}).
\end{multline*}
The last expectation is finite as we mentioned above.

\medskip

2. Let $T_r^+:=\inf\{s>0: W(s)-\ul{W}(s)\ge \max \{r,
\ol{W}(s)-\ul{W}(s)\}\}$. This is a stopping time. Introduce $B$ a
standard Brownian motion independent of $W$, and denote by $P_B$
its law. The probability in question is bounded by

\begin{multline*}
\EEE(\PPP_B(B \text{ hits first } -W^{\#}(b_r,0) \text{ and then }
3k_2 \log W^{\#}(b_r,0)|B(0)=0))=\EEE(\frac{3k_2 \log
W^{\#}(b_r,0)}{W^{\#}(b_r,0)+3k_2 \log W^{\#}(b_r,0)})\\
<\frac{6k_2}{\sqrt{r}}\EEE(W^{\#}(b_1,0)^{-1/2})<
\frac{12k_2}{\sqrt{r}}\EEE(\ol{W}(\gb_1^+)^{-1/2})
\end{multline*}

3. Let $Y:=W(b_1)-W(b_{1^++})/W^{\#}(b_1, b_{1^++})$. From Lemma
\ref{tailslemma}, $Y$ is an exponential with mean 1.

\begin{multline*} \PPP(z<m_r)\le \PPP(W(b_r)-W(b_{r^++})\le k_2\log
W^{\#}(b_r, b_{r^++})+k_2 \log r)\le \\
=\PPP(W(b_1)-W(b_{1^++})\le \frac{k_2\log (r^2 W^{\#}(b_1,
b_{1^++}))}{r}) \le \PPP(W^{\#}(b_1, b_{1^++})>r) \\ +\PPP(
W^{\#}(b_1, b_{1^++})\  Y \le \frac{k_2\log (r^3) }{r}) \le
2(2+r)e^{-r}+\PPP(Y \le \frac{3k_2\log r}{r})
<2(2+r)e^{-r}+\frac{3k_2\log r}{r}.
\end{multline*}

\end{proof}

Let $(R_i)_{i\ge1}$ be the increasing sequence of points where $b$
jumps in $[1,+\infty)$. The next lemma is the result that makes
possible to move from Lemma \ref{locprop} to the theorem.

\begin{lemma} \label{interpolation}
Let $a\in(0,1/2)$. With probability one, ultimately between any
two terms from $(R_i)_{i\ge1}$ there is at least one term from the
sequence $(\exp(k^a))_{k\ge1}$.
\end{lemma}

\begin{proof}

There are four cases for the signs of the pair $\{b_{R_i},
b_{R_{i+1}}\}$. First we show that the sequence
$(\exp(k^a))_{k\ge1}$ enters eventually in the intervals $(R_i,
R_{i+1})$ with $b_{R_i}, b_{R_{i+1}}>0$ (similarly if  $b_{R_i},
b_{R_{i+1}}<0$).

\nid For $r>0$, let
\begin{align*}
\tilde \tau_r:&=\inf\{s>0: W(s)-\ul{W}(s)=r\},\\
\tilde b_r:&=\inf\{s>0:W(s)=\ul{W}(\tilde \tau_r)\}.
\end{align*}
The process $(\tilde b_r)_{r>0}$ takes only positive values, and
it is increasing. The points where $b$ jumps from a positive to a
positive value are contained in the points where $\tilde b$ jumps.
So we will prove our claim for the process $\tilde b$. The points
where $\tilde b$ jumps in $[1,+\infty)$ make up an increasing
sequence $(h_n)_{n\ge0}$ with $h_0$ the first such point, and
$h_{n+1}:=(1+r_n)\ h_n$ for $n\ge 0$, where the $r_n$'s are i.i.d.
with density $(1+x)^{-2} \textbf{1}_{x\ge0}$ (It is the same idea
as in the proof of part 1 of Lemma \ref{thm2lemma}. It is
explained in detail in the proof of Lemma 2 of \cite{CV}). We note
that $\log (1+r_n)$ is exponential random variable with mean one.

\nid For $n\ge1$ there is a unique $k_n$ so that $c(k_n)<h_n\le
c(k_{n+1})$, i.e.,
\begin{equation} \label{interp1}
k_n^a<\log h_n\le (k_n+1)^a.
\end{equation}
We want to prove that eventually $\log h_{n+1}> (k_n+1)^a$. It is
enough to prove that $\log h_{n+1}-\log h_n> (k_n+1)^a-k_n^a$. The
last quantity is less than $a k_n^{a-1}$ (we use the fact that
$a<1$ and the mean value theorem). Also $\log h_{n+1}-\log
h_n=\log (1+r_n)$, and by the first Borel-Cantelli, we have
eventually $\log (1+r_n)>(n\log^2n)^{-1}$. So that a.s. eventually
$$
\frac{(k_n+1)^a-k_n^a} {\log h_{n+1}-\log h_n}<  ak_n^{a-1}
n\log^2n<a ((\log h_n)^{1/a}-1)^{a-1}n\log^2n.$$
In the second inequality, we used \eqref{interp1} and  $a-1<0$.
Since $\log (1+r_i)$ is exponential with mean one, we have $\log
h_n \approx n$ (for the rigorous argument we use the SLLN to say
that $\log h_n>n/2$ eventually). So that the above bound is of the
order $n^{2-1/a} \log^2 n$ which goes to zero as $n\to+\infty$
provided that $a<1/2$.

Now for the intervals $(R_i, R_{i+1})$ with $b_{R_i}
b_{R_{i+1}}<0.$

\nid For $\ell\in\D{R}$, recall the definitions of $H_\ell^-,
H_\ell^+$ given in Lemma , and moreover define $\Theta_\ell:=-\min
\{W_s: s\in[H_\ell^-, H_\ell^+]\}.$ Let $(\ell_n)_{n\ge1}$ be the
strictly increasing sequence consisting exactly of the points in
$[1,+\infty)$ where $\Theta$ jumps. At every ``time" $\ell$, we
observe $A_\ell:=W |[H_\ell^-, H_\ell^+]$. We call this a well,
and $D_\ell:=\ell+\Theta_\ell$ its depth. As $\ell$ increases, in
the picture $A_\ell$, excursions of $\overline{W}-W$ are
introduced on the right or the left. And $\Theta_l$ jumps at
$\ell$ if, just after $\ell$, an excursion is added that has
height strictly greater than $D_\ell$. For $n\ge 1$, let
$\gz_{2(n-1)}:=D_{\ell_n}, \gz_{2n-1}:=D_{\ell_n+}$. Also
$\gs_n:=\gz_{2n-1}/\gz_{2n-2}, \tau_n:=\gz_{2n}/\gz_{2n-1}.$ So
that $\gz_{2n}=\gz_0 \prod_{i=1}^n \gs_i \tau_i$ and
$\gz_{2n+1}=\gs_{n+1} \gz_{2n}$ for $n\ge0$. It can be shown that
$\{\gs_n:n\ge 1\}$ are i.i.d. with density $x^{-2}\textbf{1}_{x\ge
1}$, and $\{\tau_n:n\ge 1\}$ are i.i.d. with density $2
x^{-3}\textbf{1}_{x\ge 1}$ (see \cite{CV}, Lemma 1). If $i$ is
such that $b_{R_i} b_{R_{i+1}}<0$, then there is a $n\in \D{N}$
with $R_i<\gz_{2n+1}<\gz_{2n+2}=R_{i+1}$. As before, we prove
that, a.s. eventually, between $\gz_{2n+1}, \gz_{2n+2}$ there is a
term from the sequence $(\exp(k^a))_{k\ge1}$.
\end{proof}

\begin{lemma} \label{intervlemma}
For any $c>6$, with probability one we have
 $$(\ga_r, \gga_r)\subset(b_r-(\log b_r)^c, b_r+(\log b_r)^c)$$
 for all big $r$.
\end{lemma}
\begin{proof}

First we will show that with probability one we have
\begin{equation}\label{intervest1}(\ga_r, \gga_r)\subset(b_r-(\log r)^c, b_r+(\log
r)^c)\end{equation} for all big $r$. Define
\begin{align*}
J^+_W(r):=&\inf\{t>0:W(t)-\inf_{0\le s \le t}W(s)=r\}, \\
\gb_r^+(W):=&\inf\{t>0:W(t)=\inf_{0\le s \le J^+_W(r)}W(s)\}, \\
K^+_W(r):=&\inf\{t< \gb_r^+(W):W(t)-W(\gb_r^+(W))=r\},\\
\tilde \ga_r^+(W):=&\inf\{t>K^+_W(r):W(t)-W(\gb_r^+(W))<2k_1\log r\},\\
\tilde \gga_r^+(W):=&\sup\{t<J^+_W(r):W(t)-W(\gb_r^+(W))<2k_1\log
r\}.
\end{align*}
Then consider the process $(W^-(s))_{s\in\D{R}}$ defined by
$W^-(s)=W(-s)$ for $s\in\D{R}$, and set
\begin{align*}
\gb_r^-(W):=&-\gb_r^+(W^-), \\
\tilde \ga_r^-(W):=&-\tilde \ga_r^+(W^-),\\
\tilde \gga_r^-(W):=&-\tilde \gga_r^+(W^-).
\end{align*}
In the following we will omit $W$ in $\gb_r^+(W), \tilde
\ga_r^+(W)$, etc.

 \nid $-W(\gb_r^+), -W(\gb_r^-)$ are i.i.d with density
$r^{-1}e^{-x/r} 1_{x<0}$ (i.e., exponential with mean $r$, see
\cite{NP} Lemma of \S 1). The processes
$$
(W(K^+(r))-W(K^+(r)+t))_{0\le t\le \gb_r^+-K^+(r)},
(W(J^+(r))-W(J^+(r)-t))_{0\le t\le J^+(r)-\gb_r^+}$$
are independent three-dimensional Bessel processes starting from
zero and killed when hitting $r$.
%
%
This follows from the proof of the lemma in \S 1 of \cite{NP}, the
structure of Brownian excursions (see \cite{RY}, Chapter XII,
Theorem 4.5), and the reversibility of Brownian motion.


It holds
$$\PPP(-W(\gb_r^+)<2k_1 \log
r)< \frac{2k_1\log r}{r}$$
as $-W(\gb_r^+)$ has density bounded by $1/r$. The random
variables $\gb_r^+-\tilde \ga_r^+, \tilde \gga_r^+-\gb_r^+$ are
independent having distribution $\tau_{r}^{r-2k_1\log r}(Y)$.
Where by $\tau_c^y(Z)$ we denote the first time that a continuous
process $(Z_s)_{s\ge0}$ with $Z_0=y$ hits $c$, and $Y$ is a 3-d
Bessel process. Since $Y$ satisfies the stochastic differential
equation $dY_s=dw_s+Y_s^{-1}ds$, we have for $y, d>0$,
 $\tau_{y+d}^y(Y)\le \tau_{y+d}^y(w)$. And consequently
$$\PPP(\tau_{y+d}^y(Y)>z)\le
\PPP(\tau_{y+d}^y(w)>z)=\PPP(\tau_d^0(w)>z)=\PPP(\tau_1^0(w)>z/d^2)<d/\sqrt{z}.$$
For the last inequality, we used Fact 2. Thus for $z/\log^2 r$
large, we have
 \begin{equation}\label{locbound}\PPP(\gb_r^+-\tilde\ga_r^+>z) \le
 \PPP(\tau_{2k_1\log r}^0(w)>z) < \frac{2k_1\log r}{\sqrt{z}}.
\end{equation}
 Take $r=r_n=\exp(n^a), z=(\log r_n)^c/2$ where $c>2$. Then $z/\log^2 r=(\log
r_n)^{c-2}/2=n^{a(c-2)}/2$ is large for $n$ large, and we have the
bound
$$\PPP(\gb_{r_n}^+-\tilde \ga_{r_n}^+>\frac{(\log r_n)^c}{2})\le n^{-a(c-2)/2} 4 k_1.$$

Of course, $\PPP(\tilde\gga_{r_n}^+-\gb_{r_n}^+>(\log r_n)^c/2)$
has the same bound. Now for any $c>6$, there is an $a\in(0,1/2)$
with $-a(c-2)/2<-1$. For this choice of $a$, it holds
$\sum_{n=0}^{+\infty}n^{-a(c-2)/2}<+\infty$. Thus a.s. eventually
we have $\gb_{r_n}^+-\tilde\ga_{r_n}^+<(\log r_n)^c/2,
\tilde\gga_{r_n}^+-\gb_{r_n}^+<(\log r_n)^c/2$.

Now take an $r>0$ large. There is a unique $n$ so that $r_n<r\le
r_{n+1}$. Then $\gb_r^+=\gb^+_{r_n}$ or $\gb_r^+=\gb^+_{r_{n+1}}$
because in the interval $(r_n, r_{n+1}]$ there is at most one jump
for $\gb^+_r$ (this is included in the proof of Lemma
\ref{interpolation}). If $\gb_r^+=\gb_{r_n}^+$, then $k_1\log
r<k_1\log r_{n+1}=k_1(n+1)^a=(1+n^{-1})^ak_1\log r_n<2 k_1\log
r_n$. So $\gb_r^+-\ga_r^+<\gb_{r_n}^+-\tilde\ga_{r_n}^+<(\log
r_n)^c<(\log r)^c$, and similarly for $\gga_r^+-\gb_r^+$.

\nid If $\gb_r^+=\gb_{r_{n+1}}^+$, then
\begin{multline*}
\gb_r^+-\ga_r^+\le
\gb_{r_{n+1}}^+-\ga_{r_{n+1}}^+<\gb_{r_{n+1}}^+-\tilde\ga_{r_{n+1}}^+
<\frac{1}{2}(\log r_{n+1})^c =\frac{(n+1)^{ac}}{2}\\
=\frac{(1+n^{-1})^{ac}}{2}(\log
r_n)^c<\frac{(1+n^{-1})^{ac}}{2}(\log r)^c<(\log r)^c
\end{multline*}
 for large
$n$. And similarly for $\gga_r^+-\gb_r^+$. We do the same on the
negative side with $\ga_r^--\gb_r^-$, $\gb_r^--\gga_r^-$.

Now we claim that a.s. for all big $r$ we have
$|W(\gb_r^+)-W(\gb_r^-)|\ge2k_1 \log$. Indeed
$$\PPP(|W(\gb_r^+)-W(\gb_r^-)|<4k_1 \log
r)=\PPP(|\frac{W(\gb_r^+)}{r}-\frac{W(\gb_r^-)}{r}|<\frac{4k_1
\log r}{r})\le  \frac{8k_1\log r}{r}$$ because $W(\gb_r^+)/r,
W(\gb_r^-)/r$ are i.i.d. with density $(1-x)^{-2} 1_{x<0}$ bounded
above by 1. And since $\sum_{n=0}^{+\infty} (\log
r_n)/r_n<+\infty$, it follows that a.s. for big $n$ we have
$|W(\gb_{r_n}^+)-W(\gb_{r_n}^-)|>4k_1 \log r_n$. Similarly we show
that a.s. for big $n$ it holds
$|W(\gb_{r_n}^+)-W(\gb_{r_{n+1}}^-)|>4k_1 \log r_n$,
$|W(\gb_{r_{n+1}}^+)-W(\gb_{r_n}^-)|>4k_1 \log r_n$. And with
similar arguments as above, we show that
$|W(\gb_r^+)-W(\gb_r^-)|>2k_1 \log r$ for all big $r$, a.s.
Consequently, for all big $r$, either $(\ga_r,b_r,
\gga_r)=(\ga_r^+,\gb_r^+, \gga_r^+)$ or $(\ga_r,b_r,
\gga_r)=(\ga_r^-,\gb_r^-, \gga_r^-)$. And \eqref{intervest1}
follows from what we proved above.

To finish the proof of the lemma, it is enough to show that with
probability 1, it holds
\begin{equation}\label{lowerb}
\log |b_r|>\log r \text{ for big } r.
\end{equation}
Since for all $r$, $b_r=\gb_r^+$ or $\gb_r^-$, we will show this
for $\gb_r^+$. First we claim the following.

\nid \textsc{Claim:} There is a constant $C$ so that
$\PPP(\gb_1^+<x)\le C \sqrt{x}$ for all $x>0$.

The claim needs a proof only for small $x$. The Laplace transform
of $\gb_1^+$ is $(\sqrt{2\gl} \coth \sqrt{2\gl})^{-1}$, i.e., of
the form $\gl^{-1/2}L(\gl)$ with $L$ slowly varying function at
$+\infty$ (see \cite{NP}, Lemma of \S 1). By a Tauberian theorem
(Theorem 3 of \S XIII.5 in \cite{F}) it follows that
$\PPP(b_1^+<x)\sim x^{1/2}L(1/x)/\Gamma(3/2)$ for small $x$.

\nid Now to show the analog of \eqref{lowerb} for $\gb_r^+$, it is
enough to show that, with probability 1,
$$\gb_r^+>\frac{r^2}{\log^4 r} \text{ for all big } r.$$
We will use again an interpolation argument. This time, the
sequence $r_n=e^n$ for $n\ge1$ is enough. Observe that, because of
the above claim and scaling,
$$ \PPP(\gb_{r_n}^+<\frac{r_n^2}{\log^3 r_n})=\PPP(\gb_1^+<\frac{1}{\log^3 r_n})<C n^{-3/2}.$$
The first Borel-Cantelli lemma implies that, with probability one,
$\gb_{r_n}^+>r_n^2/\log^3 r_n$ for all big $n$. Now for $r>e$
there is unique $n$ such that $r_n<r\le r_{n+1}$, and since
$\gb^+_r\ge \gb_{r_n}^+$, we get
$$\gb_r^+ \frac{\log^4 r}{r^2}\ge \gb_{r_n}^+ \frac{\log^3 r_n}{r_n^2}
\frac{r^2_n}{r^2} \log r\ge \gb_{r_n}^+ \frac{\log^3 r_n}{r_n^2}
\frac{r^2_n}{r_{n+1}^2} \log r=\gb_{r_n}^+ \frac{\log^3
r_n}{r_n^2} \frac{\log r}{e^2}.$$ With probability one, the last
quantity is greater than one for big $r$.
\end{proof}

\bibliographystyle{annals}

\bibliography{bibliography}

\end{document}

%% file: FavP2.pstex_t
\begin{picture}(0,0)%
\includegraphics{FavP2.pstex}%
\end{picture}%
\setlength{\unitlength}{3947sp}%
\begingroup\makeatletter\ifx\SetFigFont\undefined%
\gdef\SetFigFont#1#2#3#4#5{%
  \reset@font\fontsize{#1}{#2pt}%
  \fontfamily{#3}\fontseries{#4}\fontshape{#5}%
  \selectfont}%
\fi\endgroup%
\begin{picture}(6349,3999)(289,-3823)
\put(1126,-1386){\makebox(0,0)[lb]{\smash{\SetFigFont{14}{16.8}{\rmdefault}{\mddefault}{\updefault}{\color[rgb]{0,0,0}$b_{r^-}$}%
}}}
\put(1951,-3736){\makebox(0,0)[lb]{\smash{\SetFigFont{14}{16.8}{\rmdefault}{\mddefault}{\updefault}{\color[rgb]{0,0,0}$W^{\#}(b_r, b_{r^++})$}%
}}}
\put(2351,-1386){\makebox(0,0)[lb]{\smash{\SetFigFont{14}{16.8}{\rmdefault}{\mddefault}{\updefault}{\color[rgb]{0,0,0}$b_r$}%
}}}
\put(4301,-1386){\makebox(0,0)[lb]{\smash{\SetFigFont{14}{16.8}{\rmdefault}{\mddefault}{\updefault}{\color[rgb]{0,0,0}$b_{r^++}$}%
}}}
\put(5651,-1386){\makebox(0,0)[lb]{\smash{\SetFigFont{14}{16.8}{\rmdefault}{\mddefault}{\updefault}{\color[rgb]{0,0,0}$\eta_r$}%
}}}
\put(2971,-1386){\makebox(0,0)[lb]{\smash{\SetFigFont{14}{16.8}{\rmdefault}{\mddefault}{\updefault}{\color[rgb]{0,0,0}$\zeta_r$}%
}}}
\put(721,-2461){\makebox(0,0)[lb]{\smash{\SetFigFont{14}{16.8}{\rmdefault}{\mddefault}{\updefault}{\color[rgb]{0,0,0}$W^{\#}(b_{r^-}, b_r)$}%
}}}
\end{picture}

%% file: FavP3.pstex_t
\begin{picture}(0,0)%
\includegraphics{FavP3.pstex}%
\end{picture}%
\setlength{\unitlength}{3947sp}%
\begingroup\makeatletter\ifx\SetFigFont\undefined%
\gdef\SetFigFont#1#2#3#4#5{%
  \reset@font\fontsize{#1}{#2pt}%
  \fontfamily{#3}\fontseries{#4}\fontshape{#5}%
  \selectfont}%
\fi\endgroup%
\begin{picture}(6124,4274)(139,-4023)
\put(4346,-1386){\makebox(0,0)[lb]{\smash{\SetFigFont{14}{16.8}{\rmdefault}{\mddefault}{\updefault}{\color[rgb]{0,0,0}$b_r$}%
}}}
\put(3001,-1891){\makebox(0,0)[lb]{\smash{\SetFigFont{14}{16.8}{\rmdefault}{\mddefault}{\updefault}{\color[rgb]{0,0,0}$\tilde \zeta_r$}%
}}}
\put(5721,-1891){\makebox(0,0)[lb]{\smash{\SetFigFont{14}{16.8}{\rmdefault}{\mddefault}{\updefault}{\color[rgb]{0,0,0}$\hat \zeta_r$}%
}}}
\put(4871,-1386){\makebox(0,0)[lb]{\smash{\SetFigFont{14}{16.8}{\rmdefault}{\mddefault}{\updefault}{\color[rgb]{0,0,0}$\zeta_r$}%
}}}
\put(976,-1386){\makebox(0,0)[lb]{\smash{\SetFigFont{14}{16.8}{\rmdefault}{\mddefault}{\updefault}{\color[rgb]{0,0,0}$\eta_r$}%
}}}
\put(1726,-1386){\makebox(0,0)[lb]{\smash{\SetFigFont{14}{16.8}{\rmdefault}{\mddefault}{\updefault}{\color[rgb]{0,0,0}$b_{r^++}$}%
}}}
\end{picture}

%% file: FavP4.pstex_t
\begin{picture}(0,0)%
\includegraphics{FavP4.pstex}%
\end{picture}%
\setlength{\unitlength}{3947sp}%
\begingroup\makeatletter\ifx\SetFigFont\undefined%
\gdef\SetFigFont#1#2#3#4#5{%
  \reset@font\fontsize{#1}{#2pt}%
  \fontfamily{#3}\fontseries{#4}\fontshape{#5}%
  \selectfont}%
\fi\endgroup%
\begin{picture}(5349,3999)(1414,-3748)
\put(4691,-1386){\makebox(0,0)[lb]{\smash{\SetFigFont{14}{16.8}{\rmdefault}{\mddefault}{\updefault}{\color[rgb]{0,0,0}$b_{r^++}$}%
}}}
\put(3421,-1861){\makebox(0,0)[lb]{\smash{\SetFigFont{14}{16.8}{\rmdefault}{\mddefault}{\updefault}{\color[rgb]{0,0,0}$m_r$}%
}}}
\put(4071,-1386){\makebox(0,0)[lb]{\smash{\SetFigFont{14}{16.8}{\rmdefault}{\mddefault}{\updefault}{\color[rgb]{0,0,0}$z$}%
}}}
\put(5601,-1386){\makebox(0,0)[lb]{\smash{\SetFigFont{14}{16.8}{\rmdefault}{\mddefault}{\updefault}{\color[rgb]{0,0,0}$\zeta_{r^++}$}%
}}}
\put(5226,-1386){\makebox(0,0)[lb]{\smash{\SetFigFont{14}{16.8}{\rmdefault}{\mddefault}{\updefault}{\color[rgb]{0,0,0}$\eta_r$}%
}}}
\put(2851,-1386){\makebox(0,0)[lb]{\smash{\SetFigFont{14}{16.8}{\rmdefault}{\mddefault}{\updefault}{\color[rgb]{0,0,0}$b_r$}%
}}}
\end{picture}